\def\be{\begin{equation}}
\def\ee{\end{equation}}

\def\bea{\begin{eqnarray}}
\def\eea{\end{eqnarray}}
\def\ben{\begin{eqnarray*}}
\def\een{\end{eqnarray*}}

\def\V#1{\buildrel #1\over V}
\def\L#1{\buildrel #1\over L}
\def\P#1{\buildrel #1\over P}
\def\S#1{\buildrel #1\over S}
\def\sig#1{\buildrel #1\over \sigma}

\def\Li{{\cal L}}
\def\Ve{{\cal V}}

\def\Se{{\cal S}}

\def\End{\rm End}

\def\ua{\buildrel\alpha \over u}
\def\ub{\buildrel\beta \over u}
\def\uc{\buildrel\gamma \over u}
\def\u0{\buildrel 0\over u}
\def\u#1{\buildrel #1\over u}

\def\ga{\buildrel\alpha \over g}

\def\gf#1{\buildrel #1 \over g}

\def\mua{\buildrel\alpha \over \mu}
\def\mub{\buildrel\beta \over \mu}
\def\muc{\buildrel\gamma \over \mu}

\def\muf#1{\buildrel #1 \over \mu}

\def\Wa{\buildrel\alpha \over W}
\def\Wb{\buildrel\beta \over W}
\def\Wc{\buildrel\gamma \over W}

\def\yi{\buildrel i\over y}
\def\yj{\buildrel j\over y}
\def\zi{\buildrel i\over z}

\def\cR#1{{\cal R}^{(#1)}}

\def\cUS#1{\buildrel #1\over {\cal U}}

\def\Rab{{\buildrel \alpha\beta \over R}\!}
\def\Raa{\buildrel \alpha\alpha \over R}

\def\Rff#1#2{\buildrel #1 #2 \over R}
\def\Rtff#1#2{\buildrel #1 #2 \over {\mathaccent 20 {R}}}
\def\Rpff#1#2{\buildrel #1 #2 \over {R'}}
\def\Rmff#1#2{\buildrel #1 #2 \over {R^{-1}}}
\def\Rpmff#1#2{\buildrel #1 #2 \over {{R'}^{-1}}}
\def\Rtmff#1#2{\buildrel #1 #2 \over {\mathaccent 20 {R}}^{-1}}

\def\gua{{\buildrel{\alpha} \over U}\!}
\def\guabar{{\buildrel{\bar\alpha} \over U}\!}
\def\guf#1{{\buildrel{#1} \over U}\!}

\def\Ua#1{{\buildrel\alpha \over U}(#1)}

\def\sixj#1#2#3#4#5#6{\left\{\matrix{#1 & #2 & #3 \cr #4 & #5 & #6
\cr}\right\}_q}

\def\Instate#1#2#3#4#5#6#7#8#9{ {\cal I}(#1 {\buildrel #2 \over #3} #4
{\buildrel #5 \over #6} \cdots
#7 {\buildrel #8 \over #9})}
\def\Outstate#1#2#3#4#5#6#7#8#9{ {\cal O}(#1 {\buildrel #2 \over #3} #4
{\buildrel #5 \over #6} \cdots
#7 {\buildrel #8 \over #9})}

\def\Intri#1#2#3#4#5#6{ {\cal I}(#1 {\buildrel #2 \over #3} \cdots #4
{\buildrel #5 \over #6})}
\def\Outtri#1#2#3#4#5#6{ {\cal O}(#1 {\buildrel #2 \over #3} \cdots #4
{\buildrel #5 \over #6})}

\def\Incross#1#2#3#4#5#6#7#8{ {\cal I}(#1\cdots {\buildrel #2 \over #3}
#4 {\buildrel #5 \over #6} \cdots {\buildrel #7 \over #8})}
\def\Outcross#1#2#3#4#5#6#7#8{ {\cal O}(#1\cdots {\buildrel #2 \over #3}
#4 {\buildrel #5 \over #6} \cdots {\buildrel #7 \over #8})}

\newtheorem{definition}{Definition}
\newtheorem{proposition}{Proposition}
\newtheorem{theorem}{Theorem}
\newtheorem{lemma}{Lemma}

\def\bd{\begin{definition}}
\def\ed{\end{definition}}
\def\bp{\begin{proposition}}
\def\ep{\end{proposition}}
\def\bl{\begin{lemma}}
\def\el{\end{lemma}}
\def\proof{{\sl Proof:}\,}
\def\cqfd{$\Box$}

\documentstyle[bezier,psfig]{article}
\title{Chern Simons theory on a lattice and a new description of 3-manifolds
invariants.}
\author{E.Buffenoir\thanks{e-mail: buffenoi@orphee.polytechnique.fr},\cr \cr
 Centre de Physique Theorique Ecole Polytechnique\cr
91128 Palaiseau Cedex\cr
France}

\date{\today}
\begin{document}
\maketitle

\begin{abstract}
A new approach to the quantization of Chern-Simons theory has been developed in
recent papers
\cite{FR,BR1,BR2,AGS}. It uses a
"simulation" of the moduli space of flat connections modulo the gauge group
which reveals to be
related to a lattice gauge theory based on a quantum group. After a
generalization of the formalism of
q-deformed gauge theory to the case of root of unity, we compute explicitely
the correlation
functions associated to Wilson loops (and more generally to graphs) on a
surface with punctures, which
are the interesting quantity in the study of moduli space. We then give a new
description of
Chern-Simons three manifolds invariants based on a description in terms of the
mapping class group of
a surface. At last we introduce a three dimensional lattice gauge theory based
on a quantum group which
is a lattice regularization of Chern-Simons theory. \end{abstract}

\section{Introduction}
This paper is the third part of a study of combinatorial quantization of
Chern-Simons
theory\cite{BR1,BR2}. Several ideas developed here are in fact  products of
those introduced by
V.V.Fock and A.A.Rosly in their study of  Poisson structures on the moduli
space of flat
connections\cite{FR}. To understand the motivation of these papers we must
recall some general facts
about 3D Chern-Simons theory. Chern-Simons theory is a gauge theory in 3
dimensions defined by the
action principle  \be {\cal S}_{CS}= \frac{k}{4 \pi} tr( \int_{{\cal M}}  {\cal
A}d{\cal A}+{2 \over
3}{\cal A}^3 )  \ee where ${\cal M}$ is a 3-manifold, $k$ a positive integer
and ${\cal A}$ a
connection associated to a semisimple Lie algebra ${\cal G}.$ If we first
suppose that the  manifold
locally looks like a cylinder $\Sigma \times R$, $R$ considered as the  the
time direction, we can
consider the Chern-Simons theory in an Hamiltonian point of view. We will
denote by $A$ the two
space-components of the gauge field taken to be the dynamical variables of the
theory, and the time
component $A_0$ will become a Lagrange multiplier. With these notations the
action can be written:
\be
{\cal S}_{CS}= \frac{k}{4 \pi} tr( \int_{{\cal M}} (- A \partial_0 A + 2 A_0 F
)dt).
\ee
The first
term gives the Poisson structure:
\be
\{ A^a_i(x,y), A^b_j(x',y') \} =
\frac{-2\pi}{k}\delta^{ab}\epsilon_{ij}\delta^{(2)}((x,y),(x',y')).
\label{Poisson}
\ee
The hamiltonian is a combination of constraints and the
second term imposes as a constraint that the curvature of the connection $A$ is
zero:
\be
F=dA+A^2=0.\label{constraint}
\ee
Computing the Poisson brackets of the constraints we obtain that they are first
class:
\be
\{ F^a(x,y), F^b(x',y') \} =
\frac{2\pi}{k}f^{ab}_{c}F^c(x,y)\delta^{(2)}((x,y)(x',y')).
\ee
where the $f^{ab}_{c}$ are the structure constants of the Lie algebra ${\cal
G}.$
The constraints (\ref{constraint}) generate the infinitesimal gauge
transformations of the gauge field then the phase space of the hamiltonian
Chern-Simons theory is the space of flat connections modulo gauge
transformations.\\
Quantizing the latter Poisson structure of gauge fields in the usual way,
\be
[ A^a_i(x,y), A^b_j(x',y') ] =
\frac{-2\pi}{k}\delta^{ab}\epsilon_{ij}\delta^{(2)}((x,y)(x',y'))
\ee
we are led to work with an infinite dimensional algebra, observables becoming
functionals over the
elements of this algebra. The idea of V.V.Fock and A.A.Rosly is completely
different.Let us briefly
describe it. \\
 The space of flat connections modulo the gauge group having only a finite
number of freedom degrees,
 we reduce
the connection to live on a graph encoding the topology of the surface and
equip
this "graph connection" with a Poisson structure such that the Poisson
structure induced on the gauge
 invariant observables is compatible with that induced from
the usual symplectic structure. We obtain a "simulation" of hamiltonian
Chern-Simons theory in the
sense that the operator algebra derived from both descriptions are the same.\\

We will consider a graph dividing the surface into contractile plaquettes. This
graph will be
equiped with an additionnal structure of "ciliated fat graph", i.e. a graph
with a linear order
between  adjacent links at each vertex. Let $x,y$ be neighbour vertices of the
graph, we will denote
by $U_{[x,y]}$ the  parallel transport operator associated to the link and the
connection. The gauge
group acts in the usual way on this object: \be
U_{[x,y]} \rightarrow g_{x}U_{[x,y]}g^{-1}_{y}
\ee
We can define a Lie Poisson structure on such objects \cite{FR} owning the
following fascinating
property: \\
the space of flat graph connections modulo graph gauge group is Lie Poisson
isomorphic
to the space of flat connections modulo gauge group on the surface.\\

The problem has been reduced to a finite dimensional problem by limiting the
gauge group to act at a finite number of sites.\\
The quantization of such a Lie Poisson structure leads us to an exchange
algebra
on which acts a quantum group. The object of \cite{BR1,AGS} was to define this
algebraic structure.\\

In a second time we found a projector in this algebra imposing "a posteriori"
the flatness condition,
the result was then a two dimensional lattice gauge theory based on a quantum
group. The correlation
functions associated to gauge invariant objects in this theory being related to
expectation values in
Chern-Simons theory.\\

In our second paper \cite{BR2} we further investigated the algebra of gauge
invariant elements,
particularily the algebra associated to loops. Our aim was to describe a new
approach to knots
invariants, showing in a well defined framework the relation between
Reshetikhin-Turaev invariants and
Chern-Simons theory.\\

The aim of the following paper is double:\\
\begin{itemize}
\item to  investigate further the
computation of the correlation functions of our theory on any surface and
construct the derived
invariants associated to the mapping class group, and establish a new
description of three manifold
invariants using a generalization of our previous construction to the case of
roots of unity\\
\item
 to build up a three dimensional lattice q-gauge theory associated to
triangulations of any
3-manifolds which will describe a well defined finite path integral formula for
Chern-Simons theory
and the way to compute any correlation functions of this theory.\\ This work is
a tentativ to revisit
the work of E.Witten \cite{W1} with a well defined formalism allowing a lot of
new computations and
specially computations of invariants associated to intersecting loops in
Chern-Simons theory.
\end{itemize}

\section{Lattice gauge theory based on a quantum group}

\subsection{Quantum groups and exchange algebras associated to fat graphs}

In this chapter, after a brief summary of the results of \cite{AC} \cite{MS},
we will further
develop the notion of quantum group at root of unity in the dual version and
then, using this construction, generalize the results on gauge fields algebra
developed in \cite{BR1,BR2} in a way quite different from that described in
\cite{AGS}.
We will consider a Hopf algebra $({\cal A},m,1,\Delta,S,\epsilon).$ To simplify
we will take
${\cal A}={\cal U}_q ( sl_2 )$ with $q$ being a complex number different from
$\pm 1.$ As usual we will refer to ${\bf R}$ as the universal $R-$matrix
associated to ${\cal A}$ ( we will often write ${\bf R}=\sum_i a_i \otimes b_i$
), $u$ the element defined by $u= \sum_i S(a_i)b_i$ verifying the usual
properties and $v$ the ribbon element defined by $v^2=uS(u)$ (for details see
\cite{Dr,RT1,RT2}).\\
Depending on whether $q$ is a root of unity or not, the representation theory
of ${\cal A}$ is
completely different. We will denote by $Irr({\cal A})$ the set of equivalence
classes of finite
dimensional irreducible representations of  ${\cal A}.$ In each class
$\dot{\alpha}$ we will pick
out a representativ $\alpha$ and, like often in physics, denote equivalently by
$\alpha$ or
$V_{\alpha}$the representation space associated to $\alpha.$ We will denote by
$\bar{\alpha}$ (resp.
$\tilde{\alpha}$) the right (resp. the left ) contragredient representation
build up from the
antipode (resp. the inverse of the antipode) by $\bar{\alpha}={}^{t}\alpha
\circ S$ and $0$ the one
dimensional representation associated to $\epsilon.$ The tensor product of two
representations is
defined by the coproduct $\Delta.$ If q is a root of unity the decomposition of
the tensor product
of two irreducible representations can involve indecomposable representations
(i.e representations
which are not irreducible but cannot ye!
 t be decomposed in a direct sum of  stable ${\cal A}-$modules). We are able to
introduce
$\Psi_{\alpha \beta}^{\gamma,m}$ and $\Phi_{\gamma,m}^{\alpha \beta}$
respectively projection of the
 tensor product of $\alpha$ and $\beta$ on the m-th isotypic component $\gamma$
and the inclusion of
$\gamma$ in the tensor product $\alpha \otimes \beta.$ Using this notation we
will make one more
restriction between "physical" representations verifying $\Psi_{\alpha
\bar{\alpha}}^{0}\Phi_{0}^{\alpha \bar{\alpha}} \not= 0$ and the other
representations, we will
denote by $Phys({\cal A})$ this subset of $Irr({\cal A})$. We will introduce a
new tensor product
between elements of $Phys({\cal A})$ simply realizing a truncation of the
previous one, defined by:
\be \alpha \otimes \beta = \bigoplus_{\gamma \in Phys({\cal A})}N^{\alpha
\beta}_{\gamma} \gamma.
\ee $N$ is the fusion matrix of ${\cal A}$ and we will also use the notation
$\delta(\alpha \beta
\gamma)$ to be equal to 1 or 0 depending on whether $\gamma$ occurs or not in
the decomposition of
the tensor product $\alpha \otimes \beta.$
 If q is generic the tensor product of two elements of $Irr(A)$ can be
decomposed in a direct
 sum of elements of $Irr(A).$ Moreover all irreducible representations are
"physical", so we do not
have to change the tensor product in this case. We will associate new
projection and inclusion
operators $\psi_{\alpha \beta}^{\gamma,m}$ and $\phi_{\gamma,m}^{\alpha \beta}$
build up from the
truncated tensor product. We will also use the $6-j$ notation
$\sixj{\alpha}{\beta}{\gamma}{\delta}{\mu}{\nu}$ defined in the usual way from
projection and
inclusion operators (see \cite{KR} for definitions and properties).\\ We will
use the following
notation replacing the coproduct by a truncated coproduct for any element $\xi$
of the algebra
${\cal A}$: \be {\buildrel {\alpha \otimes \beta} \over {\xi}}= \sum_{\gamma
\in Phys({\cal A})}
\phi^{\alpha \beta}_{\gamma,m} {\buildrel \gamma \over \xi}
\psi^{\gamma,m}_{\alpha \beta}. \ee The
antipode and counity maps do not change through truncation and we will again
have: \be {\buildrel
\alpha \over S(\xi)}= {}^{t} {\buildrel \bar{\alpha} \over \xi} \ee
The first trivial properties of the projection and inclusion operators are:
\bea
&&\psi_{\alpha \beta}^{{\gamma'},m'}\phi_{\gamma,m}^{\alpha
\beta}=\delta_{m,m'}\delta_{\gamma,{\gamma'}}\delta(\alpha \beta
\gamma)id_{\gamma}\\
&&\sum_{\gamma \in Phys({\cal A}),m} \phi^{\alpha
\beta}_{\gamma,m}\psi^{\gamma,m}_{\alpha \beta}={\buildrel {\alpha \otimes
\beta} \over {\bf 1}}\\
&&\phi^{\alpha 0}_{\beta}=\phi^{0 \alpha}_{\beta}=\psi_{\alpha
0}^{\beta}=\phi_{0 \alpha }^{\beta}=\delta_{\alpha \beta} id_{\alpha}.
\eea
The essential fact is that when the truncation is not trivial (i.e in the root
of unity case) the representations $((\alpha \otimes \beta)\otimes \gamma)$ and
$(\alpha \otimes (\beta \otimes \gamma))$ are no more equal but are equivalent,
the intertwiner map between them being ${\buildrel {\alpha \beta \gamma} \over
\Theta}$ defined by:
\be
{\buildrel {\alpha \beta \gamma} \over \Theta}=\sum_{\delta,\nu,\mu \in
Phys({\cal A})}
\sixj{\gamma}{\beta}{\delta}{\alpha}{\nu}{\mu}\phi_{\delta}^{\beta
\gamma}\phi_{\nu}^{\alpha \delta}\psi^{\mu \gamma}_{\nu}\psi^{\alpha
\beta}_{\mu}
\ee
( where we have omited the multiplicities to simplify the notation, it will
often be the case in the
following). We will denote by ${\buildrel {\alpha \beta \gamma} \over
{\Theta^{-1}}}$ its quasi-inverse.

We will often use the notation ${\buildrel {\alpha \beta \gamma} \over
\Theta}_{123}=\sum_{i} {\buildrel \alpha \over {\theta_i^{(1)}}} \otimes
{\buildrel \beta \over
{\theta_i^{(2)}}} \otimes {\buildrel \gamma \over {\theta_i^{(3)}}},$ and the
coproduct notations
${\buildrel {(\alpha \otimes \beta) \gamma \delta} \over
\Theta}_{1234}=\sum_{i} {\buildrel \alpha
\over {\theta_i^{(11)}}} \otimes {\buildrel \beta \over {\theta_i^{(12)}}}
\otimes {\buildrel \gamma
\over {\theta_i^{(2)}}} \otimes {\buildrel \delta \over {\theta_i^{(3)}}}...$\\

In the case where $q$ is generic ${\buildrel {\alpha \beta \gamma} \over
\Theta}$ is simply
the identity but more generally it is possible to collect some interesting
properties in the root of
unity case. Using the pentagonal identity and other trivial identities on $6-j$
symbols we can verify:
\bea {\buildrel {\alpha \beta (\gamma \otimes \delta)} \over
\Theta}\;\;{\buildrel {(\alpha \otimes
\beta) \gamma \delta} \over \Theta}&=&({\buildrel \alpha \over {\bf
1}}\otimes{\buildrel {\beta \gamma
\delta} \over \Theta}){\buildrel {\alpha (\beta \otimes \gamma) \delta} \over
\Theta}({\buildrel
{\alpha \beta \gamma} \over \Theta}\otimes{\buildrel \delta \over {\bf
1}})\label{pentagon}\\
{\buildrel {0 \alpha \beta} \over \Theta}={\buildrel { \alpha 0 \beta} \over
\Theta}&=& {\buildrel {
\alpha \beta 0} \over \Theta}= {\buildrel {\alpha \otimes \beta} \over {\bf
1}}\label{theta0} \eea
 and
other similar identities for $\Theta^{-1}.$ Moreover we have the quasi-inverse
properties:
\be
{\buildrel {\alpha \beta \gamma} \over \Theta}{\buildrel {\alpha \beta \gamma}
\over
{\Theta^{-1}}}={\buildrel {(\alpha \otimes \beta) \otimes \gamma} \over {\bf
1}} \mbox{  and   }
{\buildrel {\alpha \beta \gamma} \over {\Theta^{-1}}}{\buildrel {\alpha \beta
\gamma} \over
{\Theta}}={\buildrel {\alpha \otimes (\beta \otimes \gamma)} \over {\bf 1}} \ee
recalling that, here, ${\buildrel {\alpha \otimes \beta} \over {\bf 1}} $ is
simply a projector.
 Let us now define
intertwiners between $\alpha\otimes\beta$ and $\beta\otimes\alpha$ using our
basic objects
$\psi,\phi$: \bea &&P_{12}\; \Rab= \sum_{\gamma \in Phys({\cal A}),m}
\lambda_{\alpha \beta \gamma}^{1 \over 2}\; \phi^{\beta
\alpha}_{\gamma,m}\psi_{ \alpha \beta}^{\gamma,m}\label{Rdef}\\
&& P_{12}\; \Rpmff{\alpha}{\beta}= \sum_{\gamma \in Phys({\cal A}),m}
\lambda_{\alpha \beta \gamma}^{-{1 \over 2}}\; \phi^{\beta
\alpha}_{\gamma,m}\psi_{ \alpha \beta}^{\gamma,m}\\
\eea
where $R'=\sigma(R)$ and $\lambda_{\alpha \beta
\gamma}=(\frac{v_{\alpha}v_{\beta}}{v_{\gamma}})$
where $v_{\alpha}$ is the Drinfeld casimir, equal to $q^{C^{(2)}_{\alpha}}$,
where
$C^{(2)}_{\alpha}$ is the quadratic Casimir. We will denote in the following
$\Rab=\sum_i {\buildrel
\alpha \over a_i} \otimes  {\buildrel \beta \over b_i}$ and $\Rab^{-1}=\sum_i
{\buildrel \alpha \over
c_i} \otimes  {\buildrel \beta \over d_i}$ and use sometimes the notation
$R^{(+)}=R$ and
$R^{(-)}=R'^{-1}$  . Using the hexagonal identities on the $6-j$ symbols it can
be shown that: \bea
&&\Rff{(\alpha \otimes \beta)}{\gamma}={\buildrel {\gamma \alpha \beta} \over
\Theta}\Rff{\alpha}{\gamma}{\buildrel {\alpha \gamma \beta} \over
{\Theta^{-1}}}\Rff{\beta}{\gamma}{\buildrel {\alpha \beta \gamma} \over
\Theta}\\
&&\Rff{\alpha}{(\beta \otimes \gamma)}={\buildrel {\beta \gamma \alpha} \over
{\Theta^{-1}}}\Rff{\alpha}{\gamma}{\buildrel {\beta \alpha \gamma} \over
{\Theta}}\Rff{\alpha}{\beta}{\buildrel {\alpha \beta \gamma} \over
{\Theta^{-1}}}\label{quasitriang}
\eea which is simply the analog of the quasitriangularity property of
$R-$matrices.\\ This matrix is
no more inversible but we have: \be \Rab \Rmff{\alpha}{\beta}=
\sigma({\buildrel {\alpha \otimes
\beta} \over {\bf 1}}) \mbox{  and  } \Rmff{\alpha}{\beta}\Rab = {\buildrel
{\alpha \otimes \beta}
\over {\bf 1}}.\label{Rinv} \ee  Let us now study the properties of the
antipodal map and develop the
analog of the ribbon properties \cite{AC}.
 We will denote by ${\buildrel \alpha \over A}$ and
${\buildrel \alpha \over B}$ the matrices  defined by $\psi^{0}_{ \bar{\alpha}
\alpha}=
< . , {\buildrel \alpha \over A} .>$
and $\phi_{0}^{\alpha \bar{\alpha}}= (\lambda \rightarrow \lambda \sum_i
{\buildrel \alpha \over B}
{\buildrel \alpha \over{e_i}} \otimes {\buildrel {\bar{\alpha}} \over {e^i}}$)
where ${\buildrel
\alpha \over {e_i}}$ (resp. ${\buildrel {\bar{\alpha}} \over {e^i}}$) is a
basis of the representation
space of $\alpha$ (resp. $\bar{\alpha}$) and $<.,.>$ is the duality bracket. To
choose the
normalisation of  $\phi$ and $\psi$s we will impose the ambiant isotopy
conditions: \be \sum_i
\theta^{(1)}_i B S(\theta^{(2)}_i) A \theta^{(3)}_i = 1 \mbox{  and  } \sum_i
S({\theta^{-1}}^{(1)}_i)
A {\theta^{-1}}^{(2)}_i  B S({\theta^{-1}}^{(3)}_i) = 1 \ee
 In order to generalize the known properties relative to the antipode, we will
also introduce some
notations which will be useful in the following: \bea &&{\buildrel {\alpha
\beta} \over G } =
\sum_{i,j} ({\buildrel {\alpha} \over {S(\theta^{-1 (12)}_i)}}\otimes
{\buildrel {\beta} \over
{S(\theta^{-1 (11)}_i)}}) ({\buildrel {\alpha} \over
{S(\theta^{(2)}_j})}\otimes {\buildrel {\beta}
\over {S(\theta^{(1)}_j})})({\buildrel \alpha \over A} \otimes {\buildrel \beta
\over A} )({\buildrel
{\alpha} \over {\theta^{(3)}_j}}\otimes {\buildrel {\beta} \over {\bf
1}})({\buildrel {\alpha} \over
{\theta^{-1 (2)}_i }}\otimes {\buildrel {\beta} \over {\theta^{-1
(3)}_i)}}\nonumber\\ &&{\buildrel
{\alpha \beta} \over D } = \sum_{i,j} ({\buildrel {\alpha \otimes \beta} \over
{\theta^{(1)}_i}})
({\buildrel {\alpha} \over {\theta^{-1 (1)}_j}}\otimes {\buildrel {\beta} \over
{\theta^{-1 (2)}_j}})(
{\buildrel \alpha \over B} \otimes {\buildrel \beta \over B} ) ({\buildrel
\alpha \over {\bf 1}}
\otimes {\buildrel {\beta} \over {S(\theta^{-1(3)}_j})})  ({\buildrel {\alpha}
\over {S(\theta^{
(3)}_i)}}\otimes{\buildrel {\beta} \over {S(\theta^{ (2)}_i })}) \nonumber\\
&&{\buildrel {\alpha
\beta} \over f }= \sum_i ({\buildrel {\alpha} \over {S(\theta^{-1
(12)}_i})}\otimes {\buildrel {\beta}
\over {S(\theta^{-1 (11)}_i})}  ) {\buildrel {\alpha \beta} \over G }\;
 {\buildrel {\alpha \otimes \beta} \over {(\theta^{-1 (2)}_i B S(\theta^{-1
(3)}_i))} }
\eea
it can be shown that the latter matrices verify:
\bea
&&{\buildrel {\alpha \beta} \over {f^{-1}} }\;{\buildrel {\alpha \beta} \over G
}= {\buildrel {\alpha \otimes \beta} \over A } \mbox{  and  }
{\buildrel {\alpha \beta} \over D }\;{\buildrel {\alpha \beta} \over f }=
{\buildrel {\alpha \otimes \beta} \over B }\label{deltaA}\\
&&\phi^{\alpha \beta}_{\gamma}={\buildrel {\alpha \beta} \over
{f^{-1}}}\;{}^{t}\psi^{\bar{\gamma}}_{\bar{\beta}\bar{\alpha}}\mbox{  and  }
\psi_{\alpha \beta}^{\gamma}=
{}^{t}\phi_{\bar{\gamma}}^{\bar{\beta}\bar{\alpha}}\;{\buildrel {\alpha \beta}
\over f }
\eea
We endly introduce the element $u$ associated to the square of the antipode,
defined by:
\be
u=\sum_{i,j} S(\theta^{-1 (2)}_i B S(\theta^{-1 (3)}_i))S(b_j)A a_j\theta^{-1
(1)}_i,
\ee
u is invertible and
\bea
1&=&u\sum_{i,j} S^{-1}(\theta^{-1 (1)}_i) S^{-1}(A d_j)c_j\theta^{-1 (2)}_i B
\theta^{-1 (3)}_i=\\
&=&S^{2}(\sum_{i,j} S^{-1}(\theta^{-1 (1)}_i) S^{-1}(A d_j)c_j\theta^{-1 (2)}_i
B \theta^{-1 (3)}_i)u
\eea
moreover we have as usual the essential property
\be
\forall \xi \in {\cal A},S^2(\xi)=u \xi u^{-1}
\ee
and the usual corollaries
\bea
&&S^{2}(u)=u\\
&&uS(u)=S(u)u \mbox{  is central}\\
&&\sum_i S(b_i) A a_i=S(A)u=S(u)u\sum_i S(c_i)Ad_i\label{contract}\\
&&\epsilon(u)=1
\eea
We will denote by $v$ the element satisfying:
\bea
&&v^2=uS(u)\\
&& S(v)=v \mbox{  and  } \epsilon(v)=1.
\eea
and by $\mu$ the element $uv^{-1}.$ Then it can be shown that:
\be
{\buildrel {\alpha \otimes \beta} \over \mu} = {\buildrel {\alpha \beta} \over
{f^{-1}}}
{\buildrel {\alpha \beta} \over {(S \otimes S)(\sigma(f))}}\;\;({\buildrel
{\alpha} \over \mu}\otimes{\buildrel {\beta} \over \mu})
\ee
using the latter notations it can be shown that $\phi_{0}^{\alpha
\bar{\alpha}}=
<S(A)\mu\;.\;,\;.>$ and $\psi^{0}_{\bar{\alpha} \alpha}=(\lambda \rightarrow
\lambda \sum_i
 {\buildrel {\bar{\alpha}} \over {e^i}} \otimes \mu^{-1}S(B){\buildrel \alpha
\over {e_i}}).$ In the
following the q-dimension of the representation $\alpha$ will be defined by
$[d_{\alpha}]=tr_{\alpha}(S(A) \mu B).$

\medskip

Now, using the latter framework we can give a well defined construction of the
quantum group in the dual version for any value of $q$. As a vector space this
algebra, called $\Gamma$, is generated by $\{ {\buildrel \alpha \over {g^i_j}}
\mbox{   with $\alpha \in Phys({\cal A})$ and $i,j=1 \cdots dim(\alpha)$} \}$
and the product is simply defined by:
\be
{\buildrel \alpha \over {g}}_1 {\buildrel \beta \over {g}}_2= \sum_{\gamma \in
Phys({\cal A})}\phi_{\gamma}^{\alpha \beta}\;\;{\buildrel \gamma \over
{g}}\;\;\psi^{\gamma}_{\alpha \beta}.
\ee
This product is not associative but verify:
\be
{\buildrel {\alpha \beta \gamma} \over \Theta}_{123}(({\buildrel \alpha \over
{g}}_1 {\buildrel \beta \over {g}}_2){\buildrel \gamma \over {g}}_3)=
({\buildrel \alpha \over {g}}_1({\buildrel \beta \over {g}}_2 {\buildrel \gamma
\over {g}}_3))
{\buildrel {\alpha \beta \gamma} \over \Theta}_{123}.
\label{quasi-associativity}
\ee
Moreover we have the exchange relation:
\be
\Rab_{12} \;\;{\buildrel \alpha \over {g}}_1 {\buildrel \beta \over {g}}_2=
{\buildrel \beta \over {g}}_2 {\buildrel \alpha \over {g}}_1 \;\Rab_{12}.
\label{exchange}
\ee
This algebra can be equiped with a coproduct and a counity:
\be
\Delta({\buildrel \alpha \over {g^i_j}})=\sum_i {\buildrel \alpha \over
{g^i_k}}\; \otimes \;{\buildrel \alpha \over {g^k_j}}
\mbox{   and   }\epsilon({\buildrel \alpha \over {g^i_j}})=\delta^i_j.
\ee
Moreover it can be shown that the antipodal map $S$ defined to be the linear
map verifying $S({\buildrel \alpha \over {g^i_j}})={\buildrel {\bar{\alpha}}
\over {g^j_i}}$ owns the properties:
\bea
&& S(g)^i_k A^k_l g^l_j=A^i_j \label{antipode1}\\
&&(S(A)\mu)^k_l g^l_j S(g)^i_k=(S(A)\mu)^i_j \label{antipode2}\\
&&S({\buildrel {\beta} \over {g_2}})S({\buildrel {\alpha} \over
{g_1}})={\buildrel {\alpha \beta} \over {f}}_{12}\;S({\buildrel {\alpha} \over
{g_1}} {\buildrel {\beta} \over {g_2}})\;{\buildrel {\alpha \beta} \over
{f^{-1}}}_{12}\\
&&S^{2}({\buildrel {\alpha} \over {g}}) = {\buildrel {\alpha} \over
{\mu}}{\buildrel {\alpha} \over {g}}{\buildrel {\alpha} \over {\mu^{-1}}}
\eea

\medskip

Our aim is now to define as in our first paper the gauge theory associated
to this gauge symmetry algebra.

Let $\Sigma$ be a compact connected  oriented surface
 with boundary $\partial\Sigma$ and let ${\cal T}$ be a triangulation of
$\Sigma.$ Let  us denote by ${\cal F}$  the
oriented faces of  ${\cal T},$  by $\Li$ the set of edges counted
with their orientation. If $l$ is an interior link, $-l$ will denote
the opposite link.
We have $\Li={\Li}^{int}\cup{\Li}^{\partial \Sigma},$ where
${\Li}^{int}, {\Li}^{\partial \Sigma}$ are respectively  the set of
interior edges and boundary edges.

Finally let us also define $\Ve$ to be the set of points (vertices)
of
this triangulation, $\Ve={\Ve}^{int}\cup {\Ve}^{\partial \Sigma}$, where
${\Ve}^{int},
{\Ve}^{\partial \Sigma}$ are  respectively
the set of interior vertices and boundary vertices.

If $l$ is an oriented link it will be convenient to write $l=xy$ where $y$ is
the departure point
of $l$ and $x$ the end point of $l.$ We will write $y=d(l)$ and $x=e(l).$

\bd[gauge symmetry algebra]
Let us define for $z\in \Ve,$ the Hopf algebra  $\Gamma_z=\Gamma\times \{z\}$
 and $\hat
\Gamma=\bigotimes_{z\in \Ve} \Gamma_{z}.$ This Hopf algebra  was called in
 \cite{BR1} ``the gauge
symmetry algebra.''
\ed
If $x$ is a vertex we shall write $\ga_x$ to denote the embedding of the
 element $\ga$ in
$\Gamma_x.$

In order to define the non commutative analogue of  algebra of gauge fields we
have to
endow the triangulation with an additional structure\cite{FR}, an order between
links incident to
a same vertex, the {\sl cilium order}.

\bd[Ciliation]
A ciliation of the triangulation is an assignment of a cilium
$c_z$ to each vertex $z$   which consists in a non zero tangent  vector at z.
The orientation of the Riemann surface
defines a canonical cyclic order of the links admitting $z$ as departure or end
point. Let $l_1, l_2$ be links incident to a common vertex $z,$
the strict partial cilium order $<_{c} $ is defined by:

 $l_1<_{c}l_2$ if $l_1\not=l_2, -l_2$ and  the unoriented links
 $c_z,l_1,l_2$ appear
in the cyclic order defined
by the orientation.
\ed
If $l_1, l_2$ are incident to a same vertex $z$ we define:

$$\epsilon(l_1,l_2)=\left\{ \begin{array}{ll}
+1 &\mbox{if}\,l_1<_{c} l_2\\
 -1& \mbox{if}\,l_2<_c l_1
\end{array}
\right. $$

\bd[Gauge fields algebra]
The algebra of {\it } gauge fields \cite{AGS}\cite{BR1} $\Lambda$ is the non
associative algebra
generated by the formal variables $\ua(l)^i_j$ with $l\in \Li,\alpha\in
Phys({\cal A}),i,j=1\cdots dim(\alpha)$ and satisfying the following relations:

\smallskip

{\bf Commutation rules}
\begin{eqnarray}
& &\Rab_{12} \ua (yx)_1  \ub (yz)_2 =  \ub (yz)_2 \ua (yx)_1\label{SS}\\
 & &\ua (xy)_1 {}(S \otimes id)(\Rab_{12}) \ub (yz)_2 =\ub (yz)_2 \ua
(xy)_1\label{ES}\\
& &\ua (xy)_1  \ub (zy)_2 (S \otimes S)(\Rab_{12}) =
 \ub (zy)_2 \ua (xy)_1 \label{EE}\\
& &\, \forall \,\,(yx), (yz) \in \Li\, x\not= z \,\,\,{\rm and}\,\,\,
xy<_{c}yz\nonumber\\
& &\ua(l){\buildrel \alpha \over A}\ua(-l)={\buildrel \alpha \over
B}\label{ES=1}\\
& &\forall \,\,l \in \Li^i, \nonumber\\
& &\ua (xy)_1 \ub (zt)_2 = \ub (zt)_2  \ua (xy)_1 \label{Udisjoint}\\
& &\forall\,\, x, y, z, t \mbox{ pairwise distinct in}\, \Ve\nonumber
\end{eqnarray}

{\bf Decomposition rule}
\bea
\ua(l)_1 \ub(l)_2&=&\sum_{\gamma,m}\phi_{\gamma,m}^{\alpha,\beta}
\uc(l)\psi_{\beta,\alpha}^{\gamma,m}{\buildrel {\alpha \beta} \over {f^{-1}}}
  P_{12}, \label{UCG}\\
\u0(l)&=&1,\,\,\forall l\in \Li.
\eea

{\bf Quasi-associativity}
Let ${\cal M}_P$ be a monomial of gauge fields algebra elements with
a certain parenthesing $P.$ For each vertex $x$ of the triangulation
we construct a tensor product of representations of ${\cal A}$ by replacing
each ${\buildrel \alpha \over {u_l}}$ in the monomial by the vector space
$\alpha$ (resp. $\bar{\alpha}$  resp. $0$) depending on whether $x$ is the
endpoint (resp. departure point  resp. not element) of the edge $l,$ and
keeping the previous parenthesing. Let us consider  two different parenthesing
$P_1$ and $P_2$ of the same monomial. We can construct for both, as described
before, the corresponding vector spaces for each $x$ and deduce the intertwiner
$\Theta_x$ relating them. The relation of quasi-associativity is then simply:
\be
(\prod_{x \in {\cal V}} \Theta_x ){\cal M}_{P_1}={\cal M}_{P_2}\label{assocalg}
\ee
\ed

\bp[ Gauge covariance ] $\Lambda$ is
a right $\hat \Gamma$ comodule defined by the morphism of algebra
$\Omega:\Lambda\rightarrow \Lambda\otimes\hat\Gamma$ :
\be
\Omega(\ua(xy))=\ga_x \ua(xy) S(\ga_y).
\ee
The definition relations of the gauge fields algebra are compatible
to the coaction of the gauge symmetry algebra.
\ep

The subalgebra of gauge invariant elements of $\Lambda$ is denoted
 $\Lambda^{inv}.$

\subsection{Invariant measure, holonomies, zero-curvature projector}

It was shown
(provided some assumption on the existence of a basis of $\Lambda$ of a special
type) \cite{AGS}\cite{BR1}
that there exists a unique non zero linear form $h\in \Lambda^{\star}$
satisfying:

\begin{enumerate}
\item (invariance) $(h\otimes id)\Omega(A)= h(A)\otimes 1 \,\,\forall A\in
\Lambda$
\item (factorisation) $h((A)(B))=h(A)h(B)\\
\forall A\in\Lambda_X, \forall B\in\Lambda_Y,\forall X, Y\subset L,\,\,
(X\cup -X) \cap (Y\cup - Y)=\emptyset$
  \end{enumerate}
(we have  used the notation $\Lambda_{X}$ for $X\subset {\cal L}$ to denote the
subalgebra of $\Lambda$ generated as an algebra by $\ua_{l}$ with $l\in X$).

It can be evaluated on any element using the formula:
\be
h(\ua(x,y){}^i_j)=\delta_{\alpha,0}
\ee
where $0$ denotes the trivial representation of dimension $1,$ i.e $0$ is the
counit.

It was convenient to use the notation $\int d h$ instead of $h.$
We obtained the important formula:
\be
h(\ua(y,x)_1 (S \otimes id)(\Raa_{12}) v_{\alpha}^{-1}\ua(x,y)_2)=
{1\over [d_{\alpha}]}P_{12}{\buildrel \alpha \over B}_1.
\label{ortho}
\ee

A path $P$ (resp. a loop $P$) is a path (resp.a loop) in the graph attached
 to the triangulation of $\Sigma$, given by the collection of its vertices,
 it will also denote equivalently the
continuous curve (resp. loop) in $\Sigma$ defined by the links of $P.$
In this article we will denote by $P=[x_n,x_{n-1},\cdots,x_1,x_0]$ a link with
departure point $x_0$ and end point $x_n$.
Following the definition for links, the departure point of $P$ is denoted
$d(P)$ and its
endpoint $e(P).$ The set of vertices (resp. edges) of the path $P$ is
denoted by ${\cal V}(P)$ (resp. ${\cal L}(P)$ ), the cardinal of this set
is called the "length" of $P$ and will be denoted by $Length(P).$

Properties of path and loops such as self intersections, transverse
intersections will always be understood as properties satisfied by the
corresponding curves on $\Sigma.$

Let $P=[x_n,...,x_0]$ be a path, we defined
 the sign $\epsilon(x_i, P)$ to be $-1$ (resp. $1$) if
$x_{i-1}x_i<_{c}x_ix_{i+1}$ (resp. $x_ix_{i+1}<_{c}x_{i-1}x_i$).

\bd[Holonomies and Wilson loops]
If $P$ is a simple  path $P=[x_n,\cdots,x_0]$ with $x_0\not= x_n$, we  define
the holonomy along $P$ by
 \be
\ua_P=v_{\alpha}^{{1\over 2}\sum_{i=1}^{n-1}\epsilon(x_i, P)}(\ua(x_n
x_{n-1}){\buildrel \alpha \over A} \ua(x_{n-1}x_{n-2}){\buildrel \alpha \over
A} \cdots{\buildrel \alpha \over A}\ua(x_{1}x_{0})) .
\ee
When $C$ is a simple loop $C=[x_{n+1}=x_0, x_n,\cdots, x_0],$
we define  the holonomy along $C$  by
 \be
\ua_C=v_{\alpha}^{{1\over 2}(\sum_{i=1}^n \epsilon(x_i, C)-
\epsilon(x_0,C))}(\ua(x_0 x_{n}){\buildrel \alpha \over A}
\ua(x_{n}x_{n-1}){\buildrel \alpha \over A} \cdots{\buildrel \alpha \over
A}\ua(x_{1}x_{0})).
\ee
We define an element of $\Lambda$, called {\sl Wilson loop} attached to $C:$
\be
\Wa_C= tr_{\alpha}(S({\buildrel \alpha \over A})\mua \ua_C).\label{Wilsonloop1}
\ee
the interior parenthesing being irrelevant because the loop is simple
( the vector space attached to a vertex occurs, first with $\alpha$, second
with $\bar{\alpha}$ and
the other times in the trivial representation) and moreover we have the
relations (\ref{theta0}). We
will also use the notation $\Wa_C=\Wa_{[x_0,x_n\cdots,x_1]}.$ \ed

The properties shown in our first paper are easily generalized:

\bp[Properties of Wilson loops]
The element $\Wa_C$ is gauge invariant and moreover
it does not depend on the departure point of  the loop $C.$
Moreover it verifies the fusion equation:
\be
\Wa_C \Wb_C=\sum_{\gamma\in Phys(A)} N_{\alpha\beta}^{\gamma} \Wc_C
\ee
\ep

\medskip

\proof
The gauge invariance is quite obvious because of relations (\ref{antipode1}),
(\ref{antipode2}).
To show the cyclicity property we must put our Wilson loop in another form
called
"expanded form" in our first paper.
 Using relations (\ref{contract}),(\ref{Rinv}) we easily obtain:
\bea
&&\Wa(C)=v_{\alpha}^{-{1\over 2}(\sum_{x\in C}\epsilon(x,C))}
tr_{{\alpha}^{\otimes n}}( (S({\buildrel \alpha \over A})\mua){}^{\otimes
n}\prod_{i=n}^1
P_{ii-1}\times\\ &&\times(\prod_{i=n}^{1}\ua(x_{j+1} x_{j})_j (S \otimes id
)(R_{jj-1}^{(\epsilon(x_{j},C))}))
 \ua(x_1 x_0)_{1})\label{Wilsonloop2}.\nonumber
\eea
In this form the cyclicity invariance is obvious using the commutation
relations
(\ref{ES}).\\
The fusion relation is less trivial to show.
 We first show a lemma describing the decomposition rules for holonomies:
\be
({\buildrel \alpha \over u}_P)_1 ({\buildrel \beta \over u}_P)_2 = \phi^{\alpha
\beta}_{\gamma}\;\; {\buildrel \gamma \over u}_P \;\; \psi^{\gamma}_{\beta
\alpha} P_{12} {\buildrel {\alpha \beta} \over {f^{-1}}}_{21}.
\ee
Indeed, using (\ref{assocalg})(\ref{ES}), we easily obtain for a path
$P=[x,y,z]$:
\ben
&&(v_{\alpha}^{{1\over 2}\epsilon(y, P)}(\ua(x y){\buildrel \alpha \over A}
\ua(yz))_1\;(v_{\beta}^{{1\over 2}\epsilon(y, P)}(\ub(x y){\buildrel \beta
\over A} \ub(yz))_2=\\
&&=\sum_{i,j,k,l}v_{\alpha}^{{1\over 2}\epsilon(y, P)}v_{\beta}^{{1\over
2}\epsilon(y, P)}(\ua(x y)_1\ub(x y)_2)(S(\theta^{(1)}_l)A \theta^{(1)}_i b_j
\theta^{-1 (2)}_k \theta^{(31)}_l)_1\times\\
&&\times (S(\theta^{(2)}_l) S(\theta^{-1(1)}_k)S(a_j)S(\theta^{(2)}_i) A
\theta^{(3)}_i \theta^{-1 (3)}_k \theta^{(32)}_l)_2\;(\ua(y z)_1\ub(y z)_2)=\\
&&=v_{\alpha}^{{1\over 2}\epsilon(y, P)}v_{\beta}^{{1\over 2}\epsilon(y,
P)}(\ua(x y)_1\ub(x y)_2)\;
{\buildrel {\alpha \beta} \over G}_{21} \Rpmff{\alpha}{\beta}\;(\ua(y z)_1\ub(y
z)_2)\nonumber\\
\een
the last equality is obtained by using successively (\ref{quasitriang}) and
(\ref{pentagon}).
Now, using (\ref{deltaA})(\ref{Rdef}), we obtain the announced result for a two
links path.
Proceeding by induction we can prove it for any simple open path.\\

Let us now consider the loop $C$ as formed by two pieces $[xy]$ and $[yx]$, we
have, using the same
properties as before:
\ben
&&(v_{\alpha}^{{1\over 2}(\epsilon(y,
P)-\epsilon(x,C))}tr_{\alpha}(S({\buildrel \alpha \over A})\mua \ua(x y)
{\buildrel \alpha \over A} \ua(y x)))
(v_{\beta}^{{1\over 2}(\epsilon(y, P)-\epsilon(x,C))}tr_{\beta}(S({\buildrel
\beta \over A})\mub \ub(x y)  {\buildrel \beta \over A} \ub(y x)))=\\
&&=\sum_{{\buildrel {p,m,l,i} \over {q,n,j,k}}}(v_{\alpha}v_{\beta})^{{1\over
2}(\epsilon(y, P)-\epsilon(x,C))}tr_{\alpha \beta}((S(\theta^{-1 (2)}_p)
S(\theta^{(3)}_m) S(c_l) S( \theta^{-1 (2)}_i) S(A) \mu \theta^{-1 (1)}_i
\theta^{(1)}_m \theta^{-1 (11)}_p)_1 \times\\
&&\times(S(\theta^{-1 (3)}_p) S(A) \mu \theta^{-1 (3)}_i d_l \theta^{(2)}_m
\theta^{-1 (12)}_p)_2 (\ua(x y)_1 \ub(x y)_2)(S(\theta^{-1 (11)}_q)
S(\theta^{(1)}_n) S( \theta^{-1 (1)}_j)\times\\
&&\times  A  \theta^{-1 (2)}_j b_k \theta^{(3)}_n \theta^{-1 (2)}_q)_1
(S(\theta^{-1 (12)}_q) S(\theta^{(2)}_n) S(a_k) S(\theta^{-1 (3)}_j) A
\theta^{-1 (3)}_q)_2 (\ua(y x)_1 \ub(y x)_2))=\\
&&= (v_{\alpha}v_{\beta})^{{1\over 2}(\epsilon(y, P)-\epsilon(x,C))}tr_{\alpha
\beta}((S \otimes S)(G_{21}R_{12})(\mu \otimes \mu) (\ua(x y)_1 \ub(x y)_2)
(G_{21}R_{21}^{-1})\times\\
&&\times (\ua(y x)_1 \ub(y x)_2) )=\\
&&= \sum_{\gamma}N^{\alpha \beta}_{\gamma}v_{\gamma}^{{1\over 2}(\epsilon(y,
P)-\epsilon(x,C))}tr_{\gamma}(S({\buildrel \gamma \over A})\muc {\buildrel
\gamma \over u}(x y)  {\buildrel \gamma \over A} {\buildrel \gamma \over u}(y
x))
\een
which ends the proof of the theorem.
\cqfd

\medskip

\bp[commutation properties]
It can also be shown that
$[\Wa_C, \Wb_{C'}]=0 $
for all simple loops $C, C'$ without transverse intersections.
\ep

Although the structure of the algebra $\Lambda$ depends on the ciliation, it
has been
 shown in \cite{AGS} that the algebra $\Lambda^{inv}$ does not depend on it up
to isomorphism.
This is
completely consistent with the approach of V.V.Fock and A.A.Rosly: in their
work the graph needs to be
endowed with a structure of ciliated fat graph in order to put on the space of
graph connections
${\cal A}^l$ a structure of Poisson algebra compatible with the action of the
gauge group $G^{l}.$
However, as a Poisson algebra ${\cal A}^l/G^{l}$ is canonically isomorphic to
the space ${\cal M}^G$
of flat connections modulo the gauge group, the  Poisson structure of the
latter being independent of
any choice of r-matrix \cite{FR}.

\bd[zero-curvature projector]
We introduced a Boltzmann weight attached to any simple loop $C$ and defined
by:
\be
\delta_{C}=\sum_{\alpha\in Phys(A)}[d_{\alpha}] \Wa_{C}.
\ee
\ed
\bp
This element satisfies the flatness relation :
\begin{eqnarray}
\delta_{C}\ua_{C}{}^i_j&=&{\buildrel \alpha \over
B}^i_j\delta_{C}.\label{flatness}
\end{eqnarray}
moreover we have
\begin{eqnarray}
(\frac{(\delta_{C})}{\sum_{\alpha \in Phys({\cal A})}
[d_{\alpha}]^{2}})^{2}&=&(\frac{(\delta_{C})}{\sum_{\alpha \in Phys({\cal A})}
[d_{\alpha}]^{2}}).\label{flatness}
\end{eqnarray}
\ep

\proof
Using the same properties as in the computation of fusion relations, we obtain
:
\ben
&&\sum_{\alpha \in Phys({\cal A})}[d_{\alpha}](v_{\alpha}^{{1\over
2}(\epsilon(y, P)-\epsilon(x,C))}tr_{\alpha}(S({\buildrel \alpha \over A})\mua
\ua(x y)  {\buildrel \alpha \over A} \ua(y x)))\times\\
&&\times(v_{\beta}^{{1\over 2}(\epsilon(y, P)-\epsilon(x,C))}tr_{\beta}(
S({\buildrel \beta \over A})\mub \ub(x y)  {\buildrel \beta \over A} \ub(y
x)))=\\
&&=\sum_{{\buildrel {\alpha \in Phys({\cal A})} \over {p,m,l,i,q,n,j,k}}}
[d_{\alpha}](v_{\alpha}v_{\beta})^{{1\over 2}(\epsilon(y,
P)-\epsilon(x,C))}tr_{\alpha \beta}((S(\theta^{-1 (2)}_p) S(\theta^{(3)}_m)
S(c_l) S( \theta^{-1 (2)}_i) S(A) \mu \theta^{-1 (1)}_i \theta^{(1)}_m \times\\
&&\theta^{-1 (11)}_p)_1 (\theta^{-1 (3)}_i d_l \theta^{(2)}_m \theta^{-1
(12)}_p)_2 (\ua(x y)_1 \ub(x y)_2) (S(\theta^{-1 (11)}_q) S(\theta^{(1)}_n) S(
\theta^{-1 (1)}_j) A  \theta^{-1 (2)}_j b_k \theta^{(3)}_n \theta^{-1 (2)}_q)_1
\times\\
&&\times(S(\theta^{-1 (12)}_q) S(\theta^{(2)}_n) S(a_k) S(\theta^{-1 (3)}_j) A
\theta^{-1 (3)}_q)_2 (\ua(y x)_1 \ub(y x)_2)S(\theta^{-1 (3)}_p)_2)=\\
&&=\sum_{\alpha,\gamma \in Phys({\cal A})}[d_{\alpha}]
\sum_{p,m}
tr_{\alpha \beta}(( S(\theta^{(3)}_m) S(A) \mu \theta^{(2)}_m \theta^{-1
(12)}_p)_1 \times\\
&&\times( \theta^{(1)}_m \theta^{-1 (11)}_p)_2 \phi^{\beta \alpha}_{\gamma}
\; \uc(x y) A \uc(y x) v_{\gamma}^{{1\over 2}(\epsilon(y, P)-\epsilon(x,C))}
\;\psi^{\gamma}_{\beta \alpha} f^{-1}_{21} \; S(\theta^{-1 (2)}_p)_1
S(\theta^{-1 (3)}_p)_2)
\een
the last equality uses again the quasitriangularity properties.
Then we obtain, for any matrix $V$ in  $End(\alpha)$:

\ben
&&\delta_C \; tr_{\alpha}( S(A) \mua V \ua_C )=\\
&&=\sum_{\gamma \in Phys({\cal A})}\sum_{p,m}
tr_{\gamma}( \; (\sum_{\alpha \in Phys({\cal A})}[d_{\alpha}]
\psi^{\gamma}_{\beta \alpha}  f^{-1}_{21} (S(\theta^{-1 (2)}_p)
S(\theta^{(3)}_m) S(A) \mu \theta^{(2)}_m \theta^{-1 (12)}_p)_1 \times\\
&&\times( S(\theta^{-1 (3)}_p) S(A) \mu \; V \;\theta^{(1)}_m \theta^{-1
(11)}_p)_2 \phi^{\beta \alpha}_{\gamma} ) \; (\uc(x y) A \uc(y x)
v_{\gamma}^{{1\over 2}(\epsilon(y, P)-\epsilon(x,C))}) \; )
\een
Due to the intertwining properties and the normalizations mentioned before of
the $\phi,\psi$s
we can conclude (see \cite{AGS}) that it does exist some complex coefficients
$A(\alpha \beta
\gamma)$ such that: \bea
&&id_K \psi^{0}_{\beta \bar{\beta}} = \sum_{\alpha}A(\alpha \beta
\gamma) \psi^{\gamma}_{ \beta \alpha
}\psi^{\alpha}_{\bar{\beta} \gamma }\\ &&\mbox{   and   }\nonumber\\
&&\phi^{\gamma}_{\beta \alpha}=A(\alpha \beta
\gamma)\frac{[d_{\gamma}]}{[d_{\alpha}]}(\psi^{\alpha}_{\bar{\beta}\gamma}{\buildrel {\beta
\bar{\beta} \gamma} \over \Theta}\otimes id_{\beta}) (\phi^{\beta
\bar{\beta}}_{0}\otimes
id_{\gamma}) \eea
we obtain \ben
&&\sum_{p,m}
\sum_{\alpha \in Phys({\cal A})}[d_{\alpha}]
\psi^{\gamma}_{\beta \alpha}  f^{-1}_{21} (S(\theta^{-1 (2)}_p)
S(\theta^{(3)}_m) S(A)
\mu \theta^{(2)}_m \theta^{-1 (12)}_p)_1 \times\\
&&\times( S(\theta^{-1 (3)}_p) S(A) \mu \; V \;\theta^{(1)}_m \theta^{-1
(11)}_p)_2
\phi^{\beta \alpha}_{\gamma} =[d_{\gamma}]id_{\gamma} tr_{\beta}(S(A)\mu \; V
\; B)
\een
then for any matrix $V$ we have $\delta_C \; tr_{\alpha}({\buildrel \alpha
\over S(A)}
\mua {\buildrel \alpha \over V }\ua_C )=\delta_C \; tr_{\alpha}( {\buildrel
\alpha \over S(A)} \mua
{\buildrel \alpha \over V }{\buildrel \alpha \over  B})$, and the linear
independance of the
generators of our algebra ensures the final result.\\ The last formula of the
proposition is a
trivial consequence of the last result. \cqfd

We were led to define an element that we called $a_{YM}=\prod_{f\in {\cal
F}}\delta_{\partial
f}.$ This element is the non commutative
analogue of the projector on the space of flat connections.

In \cite{AGS}\cite{BR1} it was proved that $\delta_{\partial f}$ is a central
element of
 $\Lambda^{inv}$ and the algebra $\Lambda_{CS}=\Lambda^{inv}a_{YM}$ was shown
to be independant, up
to isomorphism, of the triangulation. The proof is based on the lemma of
decomposition rules of the
holonomies shown before and on the quite obvious property: let $C_1$ and $C_2$
be two simple
contractile loops which interiors are disjoint and with a segment $[x y]$ of
their boundary in
common: \be
\int dh(u(xy)) \Wa_{C_1} \Wb_{C_2} = \delta_{\alpha, \beta}\Wa(C_1 \# C_2)
\ee

 As  a result it was advocated that $\Lambda_{CS}$ is the algebra of
observables of the
Chern Simons theory on the manifold $\Sigma\times [0,1].$
This is supported by the topological invariance of $\Lambda_{CS}$ (i.e this
algebra depends only on
the topological structure of the surface $\Sigma$) and the flatness of the
connection.

Our aim is now to construct in the algebra $\Lambda_{CS}$ the observables
 associated to any link in $\Sigma\times [0,1].$

\subsection{Links, chord diagrams and quantum observables}
In the following subsection and in the chapter $3$ the computations will be
made in the case of
 $q$ generic
to simplify the notations but the generalization to $q$ root of unity
can be made exactly in the same way.\\
We will consider a compact connected surface $\Sigma$ with boundary $\partial
\Sigma$.
The boundary is a set of disjoint simple closed curves
 which are designed to be
"In" or "Out". Let us draw some oriented curves on the surface $\Sigma$
defining a link $L$ ,
 assuming that their boundary is contained in $\partial \Sigma$ and
with simple, transverse intersections,
with the specification of over- or undercrossing at each intersection. We will
also consider that
representations of the quantum group are attached to connected components of
the link.\\

The data ( surface with boundary + colored link ) will be called "striped
surface", the data of "In" (resp. "Out") boundary of $\Sigma$ and $L$ with
corresponding colors will be called the "In state" (resp. the "Out state") of
the striped surface.\\

 To describe such objects we will choose a  Morse function which gives a time
direction and the set
of "equitime planes" $({\cal P}_t )_{t_i \leq t \leq t_f }$ cutting the
surface.
An equitime plane $P_t$ divides the surface in two parts called respectively
"future" and "past". On any simple curve drawn on an equitime plane the time
direction give us
an orientation of the curve, if we impose moreover a departure point $x$ on
this curve we are able
to decide if a point $z$ is on the left (resp. on the right) of another point
$z',$ if $x, z, z'$
(resp. $x, z', z$) appear in the order given by this orientation. We will
consider the surface in a
canonical position defined by the following conditions. The intersection
between the surface and the
plane for $t<t_i$ or $t>t_f$ is empty. The "In" (resp. "out") boundary is
contained in the $t=t_i$
(resp. $t=t_f$)  plane. The intersection between the  $({\cal P}_t )
 _{t_i \leq t \leq t_f }$ and the surface is a set of disjoint simple closed
curves
 $(C^t_i)_{i=1,..., n(t)} $
(not necessary disjoint at the singular times), where $n(t)$ is the number of
connected components
 of $\Sigma \cap {\cal P}_t.$
 We will call "$\varphi^3-$diagram"  of the surface
a graph drawn on it which intersections with  $({\cal P}_t )
 _{t_i \leq t \leq t_f }$ determine
 departure points on each closed curve in these sets.
We impose that the $\varphi^3-$diagram never turn around any handle of the
surface.\\
Our aim is now to define a ciliated fat graph which will encode
the topology of the striped surface,i.e. this decomposition involves only
contractile plaquettes,
it is sufficiently fine to allow us to put the link on the graph in a generic
position and allow us
to distinguish two situations related by a Dehn twist of the surface. We then
decompose the surface
in  blocks, their number being chosen with respect to the  singularities of the
Morse function  (
considered as a function over the points of the surface and of the link). The
information contained
in the Morse function is not sufficient to deal with the problem of possible
non trivial cycles of
$L$ around handles of $\Sigma$. We will rule out this problem by adding
fictively two disjoint
$\varphi^3-$diagrams of $\Sigma$ to the link $L,$ the intersections between the
link and the
$\varphi^3-$diagrams will detect the rotation of the link around an handle of
$\Sigma$, we will
then refine the decomposition with respect to these datas. We will assume that
the singularities
 of the Morse function $f$, considered now as a function of the points of
the surface, points of the link and points of the $\varphi^3-$diagrams,
correspond to different times.
We will denote by $t_0=t_i,t_1,\cdots,t_{n-1},t_n=t_f$ the different instants
corresponding
to the singularities of the Morse function. We will consider a decomposition of
the surface
and of the link in "elementary blocks" ${\cal B}_0,\cdots,{\cal B}_n$
corresponding to the
subdivision
$[t_i,t_f] = [\tau_0,\tau_1] \cup [\tau_1,\tau_2] \cup
\cdots \cup [\tau_n,\tau_{n+1}]$ where
$\tau_0=t_0, \tau_1={1\over2}(t_0+t_1), \cdots,
\tau_n={1\over2}(t_{n-1}+t_{n}), \tau_{n+1}=t_n$
are called "cutting times".
An example of a striped surface with the block decomposition
described before is shown in the following figure:
\par
\centerline{\psfig{figure=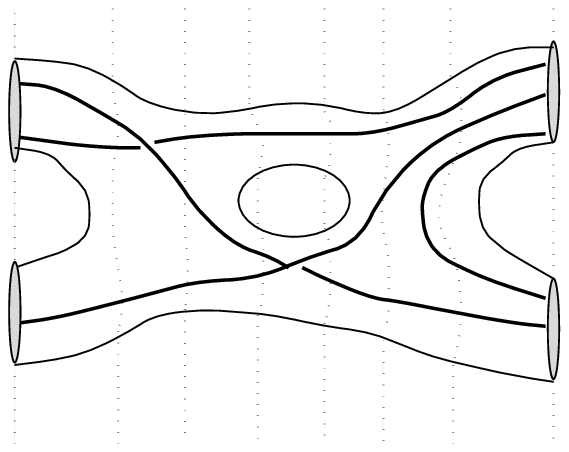}}
\par
 We will consider a triangulation ${\cal T}$ and a ciliation induced by this
block decomposition. The set ${\cal V}$ of vertices of ${\cal T}$ contains all
elements of the sets $L \cap {\cal P}_{\tau_k}$ and the singularities of the
Morse function considered as a function of the link and of the surface.
The edges of the triangulation are either a segment of the link or of the
$\varphi^3-$diagram
between two consecutiv vertices, or a segment drawn on the surface at the same
time $\tau_i$ between two emerging strings.
The plaquettes of the triangulation are the connected regions of $\Sigma$
surrounded by the edges described before.
Let us denote by $(x^k_l)$ the intersections of the link with the cutting
planes
$({\cal P}_{\tau_k}).$ The ciliation is chosen to be:
at each vertex $x^k_l$, directed to the past and to the left, just "before" the
equitime line, and
at each crossing of the oriented link, between the two outgoing strands.
We will choose a practical indexation satisfying the following properties:
$x^k_l \in C^{\tau_k}_i$ for $l \in \{ 1+\sum_{j<i}Length(C^{\tau_k}_j),
\cdots$  $\cdots, \sum_{j
\leq i}Length(C^{\tau_k}_j)\}$
 and the $x^k_l$ belonging to the same $C^{\tau_k}_i$ are ordered from the
departure point
and from left to right. The set of all $(x^k_l)$ for a given k is denoted by
${\cal
V}_{\tau_k}$. We will say that, at a vertex of $L$, the orientation of
the link is "in the sense of time" or "against the sense of time" according to
the position of
the link with respect to the
cutting plane, the two corresponding subset of ${\cal V}$ will be respectively
denoted by ${\cal V}_{-}$ and ${\cal V}_{+}.$
The set of all edges of ${\cal T}$ is again denoted by ${\cal L}$ and we will
denote by ${\cal L}_{\tau_k}$ (resp. ${\cal L}_{\tau_k < t < \tau_{k+1}}$ ) the
set of edges contained in ${\cal P}_{\tau_k}$ (resp. between ${\cal
P}_{\tau_k}$ and ${\cal P}_{\tau_k+1}$ ).
 The elements of ${\cal L}_{\tau_k < t < \tau_{k+1}}$ will be currently denoted
as $l^k_j$.
 We index the elements of this set with the following rule: begining from the
departure
point fixed by the $\varphi^3-$diagram we order the links in an obvious way
from left to right
if they do
not cross themselves and, at a crossing, the first one is that which overcross
the other.
At last we will denote by ${\cal F}_{\tau_k \le t \le \tau_{k+1}}$ the set of
faces of the triangu
lation corresponding to the
block
 ${\cal B}_k$. To each element $l$ of ${\cal L}$ we will associate a subset of
${\cal F}$ called
$Present(l)$ defined by the following rule:
if $l$ belongs to a crossing and is the overcrossing (resp.undercrossing)
strand, $Present(l)$ is
the set of the four plaquettes surrounding the crossing (resp. the set is
empty),
if $l$ is an annihilation or a creation, $Present(l)$ is the set of the two
surrounding
plaquettes,
elsewhere $Present(l)$ contains only the plaquette just at the left of $l$.
We then define $Past(l)$ to be the subset of ${\cal F}\setminus Present(l)$
such that $P \in Past(l)$ if it is on the left of $l$ in the same
block or anywhere in a past block.
A link in $\Sigma\times [0,1]$ is an embedding of $(S^1)^{\times
p'}\times([0,1])^{\times p''}$
into $\Sigma\times [0,1],$ with $\partial L \subset \partial \Sigma \times
[0,1]$.
On the set of links we can define a composition law, denoted $*$ defined as
follows:let $L$ and $L'$ be to links in $\Sigma\times[0,1]$ considered up to
ambiant isotopy. We define $L*L'$ to be the link obtained by putting $L$ in
$\Sigma \times [{1 \over 2},1]$ and $L'$ in $\Sigma \times [0,{1\over2}].$ This
composition is associative and admit the empty link as unit element. This
composition law is commutative if and only if $\Sigma$ is homeomorphic to the
sphere.
 Let us denote by $(\L{i})_{i=1\cdots p}$ the
 connected components of the link $L,$  $\alpha_i\in Phys({\cal A})$ the colour
of this component, and by $\P{i}$ the colored loop
obtained by projecting $\L{i}$ on $\Sigma.$ It is very convenient to associate
to the link $L$ a coloured chord diagram {\cal C} \cite{Va} which will encode
intersections of
the loops. This chord diagram is constructed as follows: the projection
of the link on $\Sigma$ defines $p'$ colored loops and $p''$ coloured open
paths (with boundary on
the boundary of $\Sigma$) on $\Sigma$ with transverse intersections, this
configuration of paths
defines uniquely a coloured chord diagram by the standard construction. Let us
denote by
$(\S{i})_{i=1\cdots p=p'+p''}$
the coloured circles and arcs of the chord diagram corresponding to loops and
open paths belonging
to the link ( we will call abusively "circles" the circles or the arcs of the
chord diagrams). The
family of coloured circles (resp. arcs) $(\S{i})_{i=1\cdots p'}$ (resp.
$(\S{i})_{i=p'+1\cdots p}$)
will be denoted ${\cal C}_1$ (resp. ${\cal C}_2$).
 Each circle $\S{i}$ is oriented, we
will denote by $(\yi_j)_{j=1\cdots n_i}$ the intersection points of the
circle $\S{i}$ with the chords. We will assume that they are labelled with
respect to the cyclic order
defined by the orientation of the circles.
Let $Y=\cup_{i=1}^p \{\yi_j,j=1\cdots n_i\}$, we define a relation $\sim$ on
the set $Y$
by : $y\sim y' \,\,\mbox{if and only if} \,y\,\mbox{and}\, y' \,\mbox{are
connected by a chord.}$
We will denote by $\varphi$ the
immersion of the chord diagram in $\Sigma$, in particular we have
$\P{i}=\varphi(\S{i}).$ Every intersection point of the projection of $L$ on
$\Sigma$ have exactly two inverse images by $\varphi$ in the chord diagram and
these points are  linked by a unique chord. We will denote by $\zi{}^k_j\in
\S{i}$ the
points such that  $\zi{}^k_j\in ]\yi_{j}\yi_{j-1}[$ and $\varphi(\zi{}^k_j)$ is
a vertex of the
triangulation corresponding to the cutting time $\tau_k$. We will denote by
$Z_i$ the set of all points of type $z$ in the $i-$th component, Z the union of
these sets, $Z^{\partial \Sigma}$ the subset of $Z$ formed by the points of $Z$
which belong to the boundary of $\Sigma.$
Let us denote by $\Se_i$ the family of segments forming the corresponding
circle and $\Se$ the union of these families for all components.
To each segment $s=[pq]$ we will associate two vector spaces $V_{q^{-}}$
and $V_{p^{+}}$ such that $V_{q^{-}}=
V_{p^{+}}=\V{\alpha_i}.$

$\Se$ being a finite set, let us choose on it a total ordering.
This ordering allows us to define two vector spaces $V_-$ and $V_+:$
$V_-=\bigotimes_{x\in Y\cup Z}V_{x^-}$ and
$V_+=\bigotimes_{x\in Y\cup Z}V_{x^+}$
where the order in the tensor product is taken relativ to it.

Let $a, b\in Y\cup Z$ and $\xi, \eta\in \{+, -\},$ and assume that
$\phi(a)=\phi(b),$ we will use as a shortcut the notation:
$\epsilon(a^{\xi}b^{\eta})=\epsilon(l(a^{\xi}), l(b^{\eta})).$

We define the space  $\Lambda_{\Se}$ by :
$\Lambda_{\Se}=\Lambda\otimes \bigotimes_{s\in\Se}
\End(V_{d(s)^{-}},V_{e(s)^{+}}).$

If $s$ is an element of $\Se_i$ we denote by $j_s$ the canonical injection
$j_s:\Lambda\otimes\End(V_{d(s)^{-}},V_{e(s)^{+}})\hookrightarrow
\Lambda_{\Se}.$

Let us define two types of holonomy along $s$:
$u_{s}\in\Lambda\otimes\End(V_{d(s)^{-}},V_{e(s)^{+}})$ is defined by $u_{
s}=u_{\varphi(s)},$ (the right handside has already been defined so that there
is no risk of confusion) and
$U_{s}\in\Lambda_{\Se}$ is defined by $U_{s}=j_{s}(u_{\varphi(s)}).$

\bd
Let $a, b$ two points of $Y\cup Z$ such that $\phi(a)=\phi(b)$ and define the
endomorphism:

$R^{(a^{\xi}b^{\eta})}\in \End(V_{a^{\xi}}\otimes V_{b^{\eta}})$ (resp  $
\End(V_{b^{\eta}}\otimes V_{a^{\xi}})$) if $s(a^{\xi})\triangleleft
s(b^{\eta})$ (resp if $s(b^{\eta})\triangleleft s(a^{\xi}))$ by:

$$R^{(a^{\xi}b^{\eta})}=\left\{ \begin{array}{ll}
(\alpha_1\otimes\alpha_2)(R^{(\epsilon(a^{\xi}b^{\eta}))}) &\mbox{if
$s(a^{\xi})\triangleleft s(b^{\eta})$}\\
 P_{a^{\xi}b^{\eta}}(\alpha_1\otimes
\alpha_2)(R^{(\epsilon(a^{\xi}b^{\eta}))})P_{a^{\xi}b^{\eta}}& \mbox{if
$s(b^{\eta})\triangleleft s(a^{\xi})$}
\end{array}
\right. $$
Let $<$ be any fixed strict total order on  $Y,$ we  define a family
$\{\cR{y}\}_{y\in Y}$ of elements of
$\bigotimes_{s\in\Se} End(V_{d(s)^-},V_{e(s)^{+}})$ as follows:
let $\{y, y'\}$ be any pair of points of
$Y$ such that $y\sim y',$ we can always assume (otherwise we just exchange $y$
and $y'$) that $y< y',$

\be
\cR{y}=\left\{ \begin{array}{ll}
R^{(y^{-}y^{+})-1} \,\mbox{if $\varphi (s(y))$ is above $\varphi (s(y'))$}\\
R^{(y^{+}{y'}^{+})} R^{(y^{-}y^{+})-1} R^{(y^{-}{y'}^{+})-1}\, \mbox{if
$\varphi (s(y))$ is under $\varphi (s(y'))$}
\end{array}
\right.
\ee
\be
\cR{y'}=\left\{ \begin{array}{ll}
R^{({y'}^{-}{y'}^{+})-1} \,\mbox{if $\varphi (s(y))$ is above $\varphi
(s(y'))$}\\
R^{({y'}^{-}{y'}^{+})-1}R^{({y}^{-}{y'}^{-})}R^{({y'}^{-}{y}^{+})}_{{y'}^{-}{y}^{-}}  \,  \mbox{if $\varphi (s(y))$ is under $\varphi (s(y'))$}
\end{array}
\right.
\ee
This definition defines completely the elements $\cR{y}$ for $y\in Y.$
Similarly if $z$ is an element of $Z$ we will
define $\cR{z}=R^{(z^{-}z^{+})-1}.$
\ed

We have now defined the framework necessary to associate to $L$ an element of
$\Lambda$ denoted $W_{L}$ which generalizes the construction of Wilson loops.
We denote by $<_l$ be the strict lexicographic order induced on $Y$ by the
enumeration of the connected components of $L$ and a choice of departure point
for each of these components, i.e
$\yi_{p}<_l \,\yj_{q}$  if and only if  $i<j$ or ($i=j$ and $p>q.$)

\bd
Let $P$ be a connected piece of one of the $\S{i}$. Let us choose for
simplicity  $P=[\zi_{n+1},
\yi_n, \zi_n, \cdots, \yi_1, \zi_{1}]$, we will denote the holonomy associated
to
it, by
\be
{\cal U}_P=\omega(P)U_{[\zi_{n+1}\yi_{n}]}\cR{\yi_{n}}(<_l)U_{[\yi_n \zi_n]}
\cR{\zi_n}(<_l)
\cdots \cR{\yi_1 }(<_l) U_{[\yi_1 \zi_1]},
\ee
where
$\omega(P)=v_{\alpha_i}^{-{1\over 2}
\sum_{x\in P\setminus\{\zi_{n+1}, \zi_1\}}\epsilon(\phi(x),\S{i})}.$
We will denote by $\cUS{i}$ the holonomy associated to the entire circle
$\S{i}$.
We will also define the permutation operator:
$\sigma_P = \prod_{x=\yi_n}^{\yi_1} P_{\zi_{n+1}, x}$ (where the order is given
by the order of vertices along $P$) and $\sig{i}$ will denote $\sigma_{\S{i}}.$
\ed

\bd[Generalized Holonomies and Wilson loops]
To each link in $\Sigma\times[0,1]$ we associate an element $W_L$ by the
following procedure:
let us denote by ${\cal W}_{L}$ the element
\be
{\cal W}_{L}=
\mu_{\Se}\prod_{i=1}^p\sig{i} \prod_{i=1}^{p}
\cUS{i};
\ee
where $\mu_{\Se}=\bigotimes_{x \in Z \setminus Z^{\partial \Sigma}}\mu_{x^+}$ :

 The element associated to the link $L$ is defined by
\be
W_{L}=tr_{\bigotimes_{x \in Z \setminus Z^{\partial \Sigma}}V_{x^+}}{\cal
W}_{L}
\ee
 where  $tr_{V_{+}}$ means the partial trace over the space $V_{+}$ after the
natural
identification $V_{+}=V_{-}.$
\ed

This element satisfies important properties described by the following theorem
\cite{BR2}:

\begin{theorem}
Let $L$ be a link satisfying the set of assumptions, then $W_L$ does not depend
on the labelling of the components nor does it depend on the choice of
departure points of the
 components. As a result W is a function on the space of links with values in
$\Lambda\otimes
\bigotimes_{P \in {\cal C}_2}End(V_{d(P)^-},V_{e(P)^+}).$ Moreover this mapping
is invariant under
the coaction of the gauge group at a vertex interior to the surface. If $L$ and
$L'$ are two links ,
we have the morphism property $W_{L*L'}=W_{L}W_{L'}.$ \end{theorem}

Our principal aim is the computation of the correlation function defined in an
obvious way:

\bp[Correlation functions and Ribbons invariants]
The correlation function of the link $L$ considered as immersed in
$\Sigma\times [0,1]$ is simply defined by:
\be
< W_L >_{q-YM(\Sigma)}= \int \prod_{l \in {\cal L}^{int}}dh(U_l) \;\;W_L
\;\;\prod_{F\in {\cal F}}  \delta_{\partial F}
\ee

The observable associated to $L$ will be denoted by
${\widehat W}_L=W_L\prod_{F\in {\cal F}}\delta_{\partial F}.$ This element of
$\Lambda_{CS}$ depends
only on the regular isotopy class of the link $L,$ i.e it satisfies the
Reidemeister moves of type
0,2,3. This fact was established in \cite{BR2}.\\

\par
\centerline{\psfig{figure=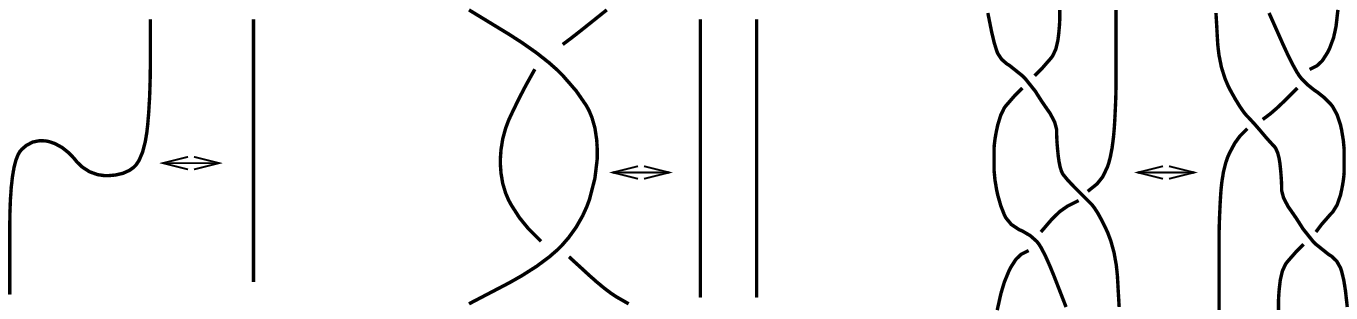}}
\par

Moreover, let  $L$ be as usual a link in $\Sigma\times [0, 1]$ and $P$ the set
of projected
curves on $\Sigma$ and  let $L^{\propto\pm}$ be another link
 whose projection $P^{\propto\pm}$ differs from $P$ by a move of type I
\par
\centerline{\psfig{figure=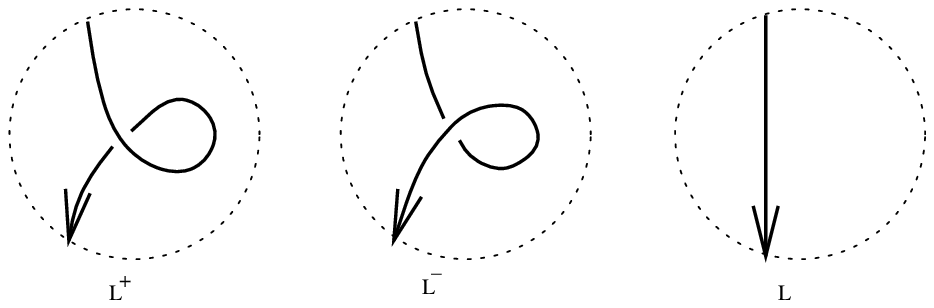}}
\par
 applied to a curve colored by $\alpha$,
we have the following relation:
\be
{\widehat W}_{L^{\propto\pm}}=v_{\alpha}^{\pm 1} {\widehat W}_{L}
\ee
The expectation value of a Wilson loop on a Riemann surface can be considered
as  an
invariant associated to a ribbon glued on the surface with the blackboard
framing.
\ep
\bigskip

\section{Computation of the correlation functions}
Our first aim is the computation of the invariants associated to links drawn on
a closed
Riemann surface embedded in $S^3.$ To realize this program we want to decompose
the computation by
introducing surfaces with boundaries and links drawn on them, already called
"striped surfaces", and
by describing the gluing operation of the latter.

This decomposition allows us to reduce the striped surface to the gluing of the
following objects,
called "elementary blocks" :  \begin{enumerate}
\item {\bf the cups}
\item {\bf the caps}
\item {\bf the (n,m)(n+m) trinions}
\item {\bf the (n+m)(n,m) trinions}
\item {\bf the free propagation of n strands}
\item {\bf the propagation of n strands with one overcrossing}
\item {\bf the propagation of n strands with one undercrossing}
 \item {\bf the (n-2)(n) creation}
\item {\bf the (n)(n-2) annihilation}
 \end{enumerate}
These objects are described in the following figure\\

\par
\centerline{\psfig{figure=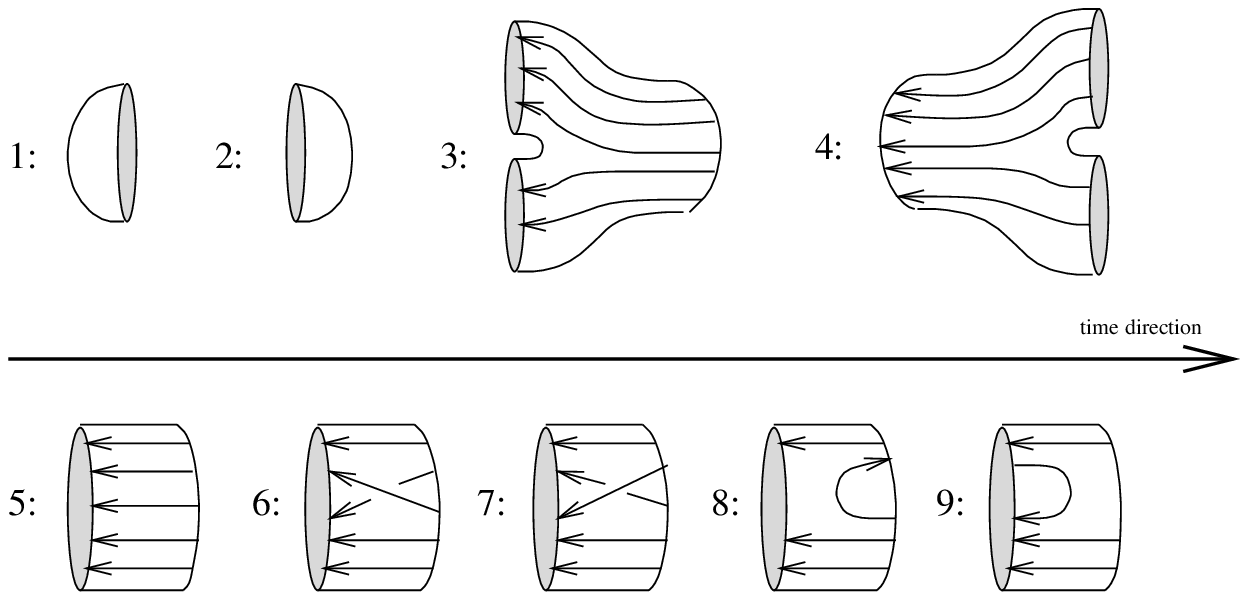}}
\par

\medskip

The correlation functions can be put in a more convenient form to reduce the
computation to the gluing of elements associated to elementary blocks.

\bl
Let us consider an element $l$ of ${\cal L}$ and an element $P$ of ${\cal F}$
then
\ben
P \in Past(l) \Rightarrow \delta_P \Ua{l} = \Ua{l} \delta_P.
\een
\el

\proof

this result is a trivial consequence of the choice of ciliation and of the
commutation properties
 developed in \cite{BR1}

\cqfd

\medskip

 From now the order induced by the orientation of the link will not be
convenient anymore, prefering
time ordering  we will introduce the vector spaces $V_{x^a}$ and $V_{x^b}$
("after" and "before")
rather than $V_{x^+}$ and $V_{x^-}.$ Let $\alpha$ denote the representation
associated to the circle
where x is taken, then $V_{x^+}=V_{x^-}=\V{\alpha}.$  If $x$ is in ${\cal
V}_+$, then
$V_{x^a}=V_{x^b}=\V{\alpha}$ and we will introduce the canonical identification
maps
$id_{(V_{x^a},V_{x^+})}$ and $id_{(V_{x^b},V_{x^-})}.$ If $x$ is in ${\cal
V}_-$, then
$V_{x^a}=V_{x^b}=\V{\bar{\alpha}}$ and we will introduce the canonical maps
${\phi^{\bar{\alpha}\alpha}_{0}}_{(V_{x^b},V_{x^+})}$ and
${\psi^{0}_{\bar{\alpha}\alpha}}_{(V_{x^a},V_{x^-})}.$
We can define a new holonomy by ${U}^{\#}_{[ x, y]}$ :
\begin{eqnarray*} {U}^{\#}_{[x, y]}&=& id_{(V_{x^a},V_{x^+})} \gua_{[ x, y ]}
id_{(V_{y^b},V_{y^-})}, \mbox{ if } x \in {\cal V}_+,  y \in {\cal V}_+\\ &=&
{\psi_{\bar{\alpha}\alpha}^{0}}_{(V_{x^b},V_{x^+})} {\gua}_{[ x, y ]}
{\phi_{0}^{\bar{\alpha}\alpha}}_{(V_{y^a},V_{y^-})}= {\guabar}_{[ y, x ]},
\mbox{ if } x \in {\cal
V}_-,  y \in {\cal V}_-\\ &=&
{\psi_{\bar{\alpha}\alpha}^{0}}_{(V_{x^b},V_{x^+})} {\gua}_{[ x, y ]}
id_{(V_{y^b},V_{y^-})}, \mbox{ if } x \in {\cal V}_-,  y \in {\cal V}_+\\ &=&
id_{(V_{x^a},V_{x^+})}
{\gua}_{[ x, y ]} {\phi_{0}^{\bar{\alpha}\alpha}}_{(V_{y^a},V_{y^-})}, \mbox{
if } x \in {\cal
V}_+,  y \in {\cal V}_- \end{eqnarray*}
 Despite its apparent complexity, this definition has a very
simple meaning. It describes the usual fact that  a strand in the direction of
the past coloured by a
representation $\alpha$  can be described by a strand in the direction of the
future coloured by a
representation $\bar{\alpha}.$

We will denote in the following:
\ben
A_k = \prod_{j=1}^{Card({\cal L}_{\tau_k \le t \le \tau_{k+1}})} ( (\prod_{P
\in Present(l^k_j)}
\delta_{\partial P})  {U}^{\#}_{l^k_j} ),
\een

the elements (  $(\prod_{P \in Present(l)}\delta_P U^{\#}_l)$ if $l$ does not
belong to a crossing
and $(\prod_{P \in Present(l)}\delta_P U^{\#}_l U^{\#}_{l'})$ if $l$ and $l'$
cross themselves )
will be called "square plaquettes" elements in the following.\\

We will also use the following permutation operator:
\ben
\sigma^{\tau_k} = \prod_{j=1}^{Card({\cal V}_{\tau_k})} (\prod_{i=Card({\cal
V}_{\tau_{k+1}})}^{1} P_{(x^k_j)^a,(x^{k+1}_i)^a} {})
\een

\bl[chronologically ordered observables]
Using these definitions, the element associated to the striped surface can be
put
in a form which respects the ordering induced by the time order:
\be
{\widehat W}_L=v_{\alpha}^{{1\over2}\sum_{x \in {\cal V}^{int}}
\epsilon(x^b,x^a)} tr_{\bigotimes_{x \in {\cal
V}^{int}}V_{x^a}}((\prod_{k=1}^{n}\sigma^{\tau_k})(\prod_{k=0}^{n}A_k))
\ee
\el
\proof
We begin with the ordering of the holonomies attached to the link.
Using the commutation relations and the properties of the $R$ matrix we obtain:
\begin{eqnarray*}
W_L&=& tr_{\bigotimes_{x \in {\cal
V}^{int}}V_{x^+}}((\prod_{k=1}^{n}\prod_{j=1}^{Card({\cal V}_{\tau_k})}
\prod_{i=Card({\cal V}_{\tau_{k+1}})}^{1} P_{(x^k_j)^+,(x^{k+1}_i)^+} )\times\\
&\times&(\otimes_{x\in {\cal V}^{int}}\mu_{x^+})
(\prod_{k=0}^{n}\prod_{j=1}^{Card({\cal L}_{\tau_k \le t \le \tau_{k+1}})}
{U}_{l^k_j}))
v_{\alpha}^{{1\over2}(\sum_{x \in {\cal V}^{int}_+} \epsilon(x^-,x^+) +
\sum_{x \in {\cal V}^{int}_-} \epsilon(x^+,x^-))}
\end{eqnarray*}
now the commutation lemma gives easily:
\ben
(\prod_{P \in {\cal F}}\delta_P)(\prod_{k=0}^{n}
\prod_{j=1}^{Card({\cal L}_{\tau_k \le t \le \tau_{k+1}})}
{U}_{l^k_j})=\prod_{k=0}^{n}\prod_{j=1}^{Card({\cal L}_{\tau_k \le t \le
\tau_{k+1}})} ( (\prod_{P
\in Present(l^k_j)} \delta_P)  {U}_{l^k_j} ), \een
and with the definition of ${U}^{\#}_{l^k_j}$
we then obtain:
\begin{eqnarray*}
{\widehat W}_L=v_{\alpha}^{{1\over2}\sum_{x \in {\cal V}^{int}}
\epsilon(x^b,x^a)}
 tr_{\bigotimes_{x \in {\cal V}^{int}}V_{x^a}}((\prod_{k=1}^{n}
\prod_{j=1}^{Card({\cal
V}_{\tau_k})} \prod_{i=Card({\cal V}_{\tau_{k+1}})}^{1}
P_{(x^k_j)^a,(x^{k+1}_i)^a} )
(\prod_{k=0}^{n}A_k))
\end{eqnarray*}
\cqfd

\bigskip

This lemma leads us to a new definition of elements associated to "striped
surfaces" which
is based on gluing chronologically ordered elementary blocks.

\bp[Correlation functions and gluing operation]
Let us consider a "striped surface" $\Sigma+L$.
Let us define an element  corresponding to  $\Sigma+L$
( which will be denoted by ${\cal A}_{\Sigma+L}$) by the following rules: \\
\begin{itemize}
\item if $\Sigma+L$ is an elementary block ${\cal B}_1$, the element of the
gauge algebra
 associated to it is:
\be
 {\cal A}_{{\cal B}_1}=\int \prod_{l \in {\cal L}_{\left] t,t' \right[ }}
dh(U_l)
\prod_{j=1}^{Card({\cal L}_1)} ( (\prod_{P \in Present(l^1_j)} \delta_P)
{U}^{\#}_{l^1_j} )
v_{\alpha}^{{1\over2}\sum_{x \in {\cal V}^{t'}} \epsilon(x^b,x^a)} \ee
\item if $\Sigma+L$ is a disjoint union of $N$ elementary blocks placed between
$t$ and $t'$ then
the element of the algebra associated to  $\Sigma+L$ is obviously the product
of the elements
associated to each elementary block, the order between them being irrelevant
because they are
commuting.\\

\item if there exists a time $t''$ between $t$ and $t'$ such that $\Sigma+L$ is
obtained by gluing
two "striped surfaces" $\Sigma_1+L_1$ and $\Sigma_2+L_2$ placed respectively
between $t$ and $t''$,
and between $t''$ and $t'.$ The element associated to $\Sigma+L$ will be
defined by: \bea
{\cal A}_{\Sigma+L}&=&{\cal A}_{\Sigma_1+L_1} \circ {\cal A}_{\Sigma_2+L_2}\\
&=& \int \prod_{l \in {\cal L}_{t''}} dh(U_l) tr_{\bigotimes_{x \in {\cal
V}^{t''}}V_{x^a}}(\sigma^{t''}{\cal A}_{\Sigma_1+L_1} {\cal
A}_{\Sigma_2+L_2})\nonumber
\eea
( the canonical choice of ciliation defined for any striped surface is
obviously compatible with the gluing operation )\\
\end{itemize}
these properties give us a new way to compute the invariants associated to
links on a closed
surface:
\be
< W_L >_{q-YM(\Sigma)}= {\cal A}_{\Sigma+L}.
\ee
\ep

\medskip

 From now the computation of the correlation functions
is reduced to the computation of elements
associated to elementary blocks. After some definitions we
will give the result of the explicit computation of these
elements.

\bd[In and Out states]

Let us consider a striped surface $\Sigma+L$.
A connected component of its "In state" is a simple loop ${\cal
C}=[x_{n+1}=x_1, x_n, \cdots, x_1]$
 oriented in the inverse clockwise sense with $n+1$ strands going through it at
each $x_i$
in the direction of the past with a representation $\alpha_i.$\\
 ( a strand in the
direction of the future with a representation $\alpha$ is reversed to the past
by changing its representation in $\bar{\alpha}$). The ciliation at each $x_i$
is chosen as in the general construction of striped surfaces.
Then, choosing $n$ other representations $(\beta_i)_{i=1,...,n}$  we define
${\cal O} \in \Lambda \otimes End(\otimes_{x \in {\cal C}} V_{x^b},{\bf C})$
and
${\cal I} \in \Lambda \otimes End({\bf C},\otimes_{x \in {\cal C}} V_{x^a})$:\\
if there is at least one emerging strand through ${\cal C}$,
\begin{eqnarray*}
\Outstate{\beta_n
}{\alpha_n}{x_n}{\beta_{n-1}}{\alpha_{n-1}}{x_{n-1}}{\beta_1}{\alpha_1}{x_1}&=&
v^{-1}_{\beta_n}tr_{\V{\beta_n}}
(\psi^{\beta_n}_{\beta_1
\alpha_1}\guf{\beta_1}_{[x_1,x_2]}\psi^{\beta_1}_{\beta_2 \alpha_2}
\cdots  \psi^{\beta_{n-1}}_{\beta_n\alpha_{n}}\guf{\beta_n}_{[x_{n},x_1]}
\Rmff{\beta_n}{\alpha_1}\muf{\beta_n})\\
\Instate{\beta_n}{\alpha_n}{x_n}{\beta_{n-1}}{\alpha_{n-1}}{x_{n-1}}{\beta_1}{\alpha_1}{x_1}&=&
v_{\beta_n}tr_{\V{\beta_n}}(\muf{\beta_n}\Rff{\beta_n}{\alpha_1}
\guf{\beta_n}_{[x_1,x_{n}]}\phi^{\beta_n \alpha_{n}}_{\beta_{n-1}} \cdots
\phi^{\beta_2
\alpha_2}_{\beta_1}\guf{\beta_1}_{[x_2,x_{1}]}\phi^{\beta_1
\alpha_1}_{\beta_n}) \end{eqnarray*}
Here and in the following we will often forget the multiplicities $m_i$ for
readibility.\\ If the
connected component of this "In state" has no emerging strand we will define
${\cal O}$ and ${\cal
I}$ to be: \be {\cal O}(\beta_0)=W^{\beta_0}_{C^{-1}}\mbox{   and   }{\cal
I}(\beta_0)=W^{\beta_0}_{C} \ee
\ed

\bigskip

The properties of the latter objects are described in the following lemma.
\bl
The properties of the In and Out states are generalizations of those of Wilson
loops.\\
Cyclicity
\bea
\Outstate{\beta_n}{\alpha_n}{x_n}{\beta_{n-1}}{\alpha_{n-1}}{x_{n-1}}{\beta_1}{\alpha_1}{x_1}&=&\Outstate{\beta_{n-1}}{\alpha_{n-1}}{x_{n-1}}{\beta_1}{\alpha_1}{x_1}{\beta_n}{\alpha_n}{x_n}\nonumber\\
\Instate{\beta_n}{\alpha_n}{x_n}{\beta_{n-1}}{\alpha_{n-1}}{x_{n-1}}{\beta_1}{\alpha_1}{x_1}&=&\Instate{\beta_{n-1}}{\alpha_{n-1}}{x_{n-1}}{\beta_1}{\alpha_1}{x_1}{\beta_n}{\alpha_n}{x_n}\nonumber\\
\eea
Gauge transformation
\bea
\Omega(\Outstate{\beta_n}{\alpha_n}{x_n}{\beta_{n-1}}{\alpha_{n-1}}{x_{n-1}}{\beta_1}{\alpha_1}{x_1})&=&\Outstate{\beta_n}{\alpha_n}{x_n}{\beta_{n-1}}{\alpha_{n-1}}{x_{n-1}}{\beta_1}{\alpha_1}{x_1} \; \prod_{i=1}^{n} S(\gf{\alpha_i}_{x_i})\nonumber\\
\Omega(\Instate{\beta_n}{\alpha_n}{x_n}{\beta_{n-1}}{\alpha_{n-1}}{x_{n-1}}{\beta_1}{\alpha_1}{x_1})&=& \prod_{i=1}^{n} \gf{\alpha_i}_{x_i}\; \Instate{\beta_n}{\alpha_n}{x_n}{\beta_{n-1}}{\alpha_{n-1}}{x_{n-1}}{\beta_1}{\alpha_1}{x_1}\nonumber\\
\eea
Scalar product
\bea
\int \prod_{l \in {\cal C}}dh(U_l)&&tr_{\otimes_{i}V_{x_i^a}}(\sigma_{{\cal C}}
\Outstate{\beta_n}{\alpha_n}{x_n}{\beta_{n-1}}{\alpha_{n-1}}{x_{n-1}}{\beta_1}{\alpha_1}{x_1}\times\\
&&\times\Instate{{\beta'}_n}{\alpha_n}{x_n}{{\beta'}_{n-1}}{\alpha_{n-1}}{x_{n-1}}{{\beta'}_1}
{\alpha_1}{x_1})\;\;=\;\;
\prod_{i=1}^{n}\delta_{\beta_i,{\beta'}_i}\nonumber \eea
\el

\proof

The cyclicity property is not completely obvious. We give here a detailed proof
of this fact:
\ben
&&\Outstate{\beta_n}{\alpha_n}{x_n}{\beta_{n-1}}{\alpha_{n-1}}{x_{n-1}}{\beta_1}{\alpha_1}{x_1}=\\
&&=v^{-1}_{\beta_n}tr_{V_{\beta_n}}
(\psi^{\beta_n}_{\beta_1
\alpha_1}\guf{\beta_1}_{[x_1,x_2]}\psi^{\beta_1}_{\beta_2
\alpha_2}\guf{\beta_2}_{[x_2,x_3]} \cdots
\psi^{\beta_{n-1}}_{\beta_n\alpha_{n}}\guf{\beta_n}_{[x_{n},x_1]}
\Rmff{\beta_n}{\alpha_1}\muf{\beta_n})\\
&&=v^{-1}_{\beta_n}tr_{V_{\beta_1}}
(\guf{\beta_1}_{[x_1,x_2]}\psi^{\beta_1}_{\beta_2 \alpha_2}
\guf{\beta_2}_{[x_2,x_3]} \cdots
\psi^{\beta_{n-1}}_{\beta_n\alpha_{n}}\guf{\beta_n}_{[x_{n},x_1]}
\psi^{\beta_n}_{\beta_1 \alpha_1}
\Rmff{\beta_1}{\alpha_1}\muf{\beta_1}v_{\alpha_2})\\
&&=v^{-1}_{\beta_n}tr_{V_{\beta_1}\otimes {V'}_{\beta_1}}
(P_{V_{\beta_1},{V'}_{\beta_1}}\guf{\beta_1}_{[x_1,x_2]}\psi^{{\beta'}_1}_{\beta_2 \alpha_2}
\guf{\beta_2}_{[x_2,x_3]}  \cdots
\psi^{\beta_{n-1}}_{\beta_n\alpha_{n}}\guf{\beta_n}_{[x_{n},x_1]}
\psi^{\beta_n}_{{\beta'}_1 \alpha_1}
\Rmff{{\beta'}_1}{\alpha_1}\muf{{\beta'}_1}v_{\alpha_2})\\
&&=\sum_{(i),(j)}v^{-1}_{\beta_n}tr_{V_{\beta_1}}
(\psi^{\beta_1}_{\beta_2 \alpha_2} b_{(i)}^{\beta_2} \guf{\beta_2}_{[x_2,x_3]}
\cdots
\psi^{\beta_{n-1}}_{\beta_n\alpha_{n}}\guf{\beta_n}_{[x_{n},x_1]}
S(a_{(j)}^{\beta_n})
\psi^{\beta_n}_{\beta_1 \alpha_1}
\Rmff{\beta_1}{\alpha_1}\muf{\beta_1}v_{\alpha_2}
b_{(j)}^{\beta_1}\guf{\beta_1}_{[x_1,x_2]} S^{2}(a_{(i)}^{\beta_1}))\\
&&=\sum_{(i)}v^{-1}_{\beta_n}tr_{V_{\beta_1}}
(S^{2}(a_{(i)}^{\beta_1})b_{(i)}^{\beta_1}
\psi^{\beta_1}_{\beta_2 \alpha_2} \guf{\beta_2}_{[x_2,x_3]}  \cdots
\psi^{\beta_{n-1}}_{\beta_n\alpha_{n}}\guf{\beta_n}_{[x_{n},x_1]}\psi^{\beta_n}_{\beta_1 \alpha_1}
\guf{\beta_1}_{[x_1,x_2]} \Rmff{\beta_1}{\alpha_2}v_{\beta_n})\\
&&=\sum_{(i)}v_{\beta_1}^{-1}tr_{V_{\beta_1}} (\psi^{\beta_1}_{\beta_2
\alpha_2}
\guf{\beta_2}_{[x_2,x_3]}  \cdots
\psi^{\beta_{n-1}}_{\beta_n\alpha_{n}}\guf{\beta_n}_{[x_{n},x_1]}\psi^{\beta_n}_{\beta_1 \alpha_1}
\guf{\beta_1}_{[x_1,x_2]} \Rmff{\beta_1}{\alpha_2}\mu_{\beta_1})\\
&&=\Outstate{\beta_1}{\alpha_1}{x_1}{\beta_{n}}{\alpha_{n}}{x_{n}}{\beta_2}{\alpha_2}{x_2} \een The
gauge transformation is very simple to derive using the decomposition rules of
the elements of the
group and we can proove the scalar product property using simply the
integration formula and the
unitarity relations of Clebsch-Gordan maps. \ben &&\int \prod_{l \in {\cal
C}}dh(U_l)tr_{\otimes_{i}V_{x_i^a}}(\sigma_{{\cal C}} \Outstate{\beta_n
m_n}{\alpha_n}{x_n}{\beta_{n-1} m_{n-1}}{\alpha_{n-1}}{x_{n-1}}{\beta_1
m_1}{\alpha_1}{x_1}\times\\
&&\;\;\;\;\;\;\;\;\;\;\;\;\;\;\Instate{{\beta'}_n
{m'}_n}{\alpha_n}{x_n}{{\beta'}_{n-1}
{m'}_{n-1}}{\alpha_{n-1}}{x_{n-1}}{{\beta'}_1 {m'}_1}{\alpha_1}{x_1})=\\
&&=\int \prod_{l \in {\cal
C}}dh(U_l) v^{-1}_{\beta_n}v_{{\beta'}_n}tr_{V_{\beta_n} \otimes
V_{{\beta'}_n}} (\psi^{\beta_n
m_n}_{\beta_1 \alpha_1}\guf{\beta_1}_{[x_1,x_2]}\psi^{\beta_1 m_1}_{\beta_2
\alpha_2} \cdots
\psi^{\beta_{n-1}
m_{n-1}}_{\beta_n\alpha_{n}}\muf{{\beta'}_n}\Rff{{\beta'}_n}{\alpha_1}
\guf{\beta_n}_{[x_{n},x_1]}\times\\ &&\times\guf{{\beta'}_n}_{[x_1,x_{n}]}
 \Rff{\beta_n}{\alpha_1}^{-1}\muf{\beta_n}
\phi^{{\beta'}_n \alpha_{n}}_{{\beta'}_{n-1} {m'}_{n-1}} \cdots
\phi^{{\beta'}_2 \alpha_2}_{{\beta'}_1
{m'}_1}\guf{{\beta'}_1}_{[x_2,x_{1}]}\phi^{{\beta'}_1
\alpha_1}_{{\beta'}_n {m'}_n})=\\ &&=\int \prod_{l \in {\cal
C}\setminus{[x_1,x_n]}}dh(U_l)
\frac{\delta_{\beta_n,{\beta'}_n}
\delta_{m_n,{m'}_n}}{[d_{\beta_n}]}tr_{V_{\beta_n}}
(\muf{\beta_n}\psi^{\beta_n m_n}_{\beta_1
\alpha_1}\guf{\beta_1}_{[x_1,x_2]}\psi^{\beta_1
m_1}_{\beta_2 \alpha_2} \cdots  \psi^{\beta_{n-1} m_{n-1}}_{\beta_n\alpha_{n}}
\phi^{{\beta}_n
\alpha_{n}}_{{\beta'}_{n-1} {m'}_{n-1}} \cdots \\ &&\cdots \phi^{{\beta'}_2
\alpha_2}_{{\beta'}_1
{m'}_1}\guf{{\beta'}_1}_{[x_2,x_{1}]}\phi^{{\beta'}_1 \alpha_1}_{{\beta'}_n
{m'}_n})=\\
&&=\frac{\delta_{\beta_n,{\beta'}_n}\delta_{m_n,{m'}_n}}{[d_{\beta_n}]}tr_{V_{\beta_n}}
(\muf{\beta_n})\prod_{i=1}^{n-1}(\delta_{\beta_i,{\beta'}_i}\delta_{m_i,{m'}_i})=\prod_{i=1}^{n}(\delta_{\beta_i,{\beta'}_i}\delta_{m_i,{m'}_i}).\\
\een This ends the proof of the lemma.\\

\cqfd

All elements associated to elementary blocks can be computed in terms of "In"
and "Out" states of the latter form.
\bp
Let us give here the expression of the elements associated to the elementary
blocks
enumerated before:\\
\begin{eqnarray}
&&{\bf {\cal A}^{elem}_{cup}}= \sum_{\beta_0} [d_{\beta_0}] {\cal O}(\beta_0)\\
&&{\bf {\cal A}^{elem}_{cap}}= \sum_{\beta_0} [d_{\beta_0}] {\cal I}(\beta_0)\\
&&{\bf {\cal A}^{elem}_{(n,m)(n+m) tri.}}= \sum_{\beta_1,\cdots, \beta_{n+m}}
[d_{\beta_{n+m}}]^{-1}
\Intri{{\beta'}_n}{{\alpha'}_n}{{x'}_n}{{\beta'}_1}{{\alpha'}_1}{{x'}_1}\Intri{{\beta''}_n}{{\alpha''}_n}{{x''}_n}{{\beta''}_1}{{\alpha''}_1}{{x''}_1}
\times\\
&&\times\Outtri{{\beta}_n}{{\alpha}_n}{{x}_n}{{\beta}_1}{{\alpha}_1}{{x}_1}
\delta_{{\beta'}_n, {\beta''}_m, {\beta}_{n+m}, \beta_m}
\prod_{k=1}^{n}\delta_{\alpha_{m+k}, {\alpha'}_{k}}
\prod_{k=1}^{m}\delta_{\alpha_{k}, {\alpha''}_{k}}
\prod_{k=1}^{n-1}\delta_{\alpha_{m+k}, {\alpha'}_{k}}
\prod_{k=1}^{m-1}\delta_{\alpha_{k}, {\alpha''}_{k}}\nonumber\\
&&{\bf {\cal A}^{elem}_{(n+m)(n,m) tri.}}= \sum_{\beta_1,\cdots, \beta_{n+m}}
[d_{\beta_{n+m}}]^{-1}
\Intri{{\beta}_n}{{\alpha}_n}{{x}_n}{{\beta}_1}{{\alpha}_1}{{x}_1}
\Outtri{{\beta'}_n}{{\alpha'}_n}{{x'}_n}{{\beta'}_1}{{\alpha'}_1}{{x'}_1}\times\\
&&\times
\Outtri{{\beta''}_n}{{\alpha''}_n}{{x''}_n}{{\beta''}_1}{{\alpha''}_1}{{x''}_1}
\delta_{{\beta'}_n, {\beta''}_m, {\beta}_{n+m}, \beta_m}
\prod_{k=1}^{n}\delta_{\alpha_{m+k}, {\alpha'}_{k}}
\prod_{k=1}^{m}\delta_{\alpha_{k}, {\alpha''}_{k}}
\prod_{k=1}^{n-1}\delta_{\alpha_{m+k}, {\alpha'}_{k}}
\prod_{k=1}^{m-1}\delta_{\alpha_{k}, {\alpha''}_{k}}\nonumber\\
&&{\bf {\cal A}^{elem}_{free}}=\sum_{\beta_1,\cdots,\beta_n}
\Intri{\beta_n}{\alpha_{n}}{x_{n}}{\beta_1}{\alpha_1}{x_1}
\Outtri{{\beta'}_n}{{\alpha'}_{n}}{{x'}_{n}}{{\beta'}_1}{{\alpha'}_1}{{x'}_1}
\prod_{i=1}^{n}\delta_{\beta_i,{\beta'}_i}\prod_{i=1}^{n}\delta_{\alpha_i,{\alpha'}_i}\\
&&{\bf {\cal A}^{elem}_{creation}}=\sum_{{\beta'}_1,\cdots,{\beta'}_n}
\Incross{\beta_n}{\alpha_{k+2}}{x_{k+2}}{\beta_k}{\alpha_{k-1}}{x_{k-1}}{\alpha_1}{x_1}
\Outtri{{\beta'}_n}{{\alpha'}_{n}}{{x'}_{n}}{{\beta'}_1}{{\alpha'}_1}{{x'}_1}\times\\
&&\times\prod_{i=k+2}^{n}(\delta_{\beta_i,{\beta'}_i}\delta_{\alpha_i,{\alpha'}_i})
\delta_{{\beta'}_{k+1},\beta_k,{\beta'}_{k-1}}\delta_{{\alpha'}_{k+1},{\bar{\alpha'}}_{k}}
(v_{{\beta'}_k} v_{{\beta}_k})^{{1\over 2}}
(\frac{[d_{{\beta'}_k}]}{[d_{{\beta}_k}]})^{{1\over 2}}
N^{{\beta'}_k,{m'}_k}_{{\beta'}_{k-1}\alpha_k} \prod_{i=1}^{k-1}
\delta_{\alpha_i,{\alpha'}_i}
\prod_{i=1}^{k-2}  \delta_{\beta_i,{\beta'}_i}\nonumber\\ &&{\bf {\cal
A}^{elem}_{annihil.}}
=\sum_{{\beta'}_1,\cdots,{\beta'}_n}
\Intri{{\beta'}_n}{{\alpha'}_{n}}{{x'}_{n}}{{\beta'}_1}{{\alpha'}_1}{{x'}_1}
\Outcross{\beta_n}{\alpha_{k+2}}{x_{k+2}}{\beta_k}{\alpha_{k-1}}{x_{k-1}}{\alpha_1}{x_1} \times\\
&&\times\prod_{i=k+2}^{n}(\delta_{\beta_i,{\beta'}_i}\delta_{\alpha_i,{\alpha'}_i})
\delta_{{\beta'}_{k+1},\beta_k,{\beta'}_{k-1}}\delta_{{\alpha'}_{k+1},{\bar{\alpha'}}_{k}}
(v_{{\beta'}_k} v_{{\beta}_k})^{-{1\over 2}}
(\frac{[d_{{\beta'}_k}]}{[d_{{\beta}_k}]})^{{1\over 2}}
 N^{{\beta'}_k,{m'}_k}_{{\beta'}_{k-1}\alpha_k} \prod_{i=1}^{k-1}
\delta_{\alpha_i,{\alpha'}_i} \prod_{i=1}^{k-2}
\delta_{\beta_i,{\beta'}_i}\nonumber\\ &&{\bf {\cal
A}^{elem}_{overcross.}}=\sum_{\beta_1,\cdots,\beta_n,{\beta'}_k}
\Incross{\beta_n}{\alpha_{k+1}}{x_{k+1}}{\beta_k}{\alpha_k}{x_k}{\alpha_1}{x_1}
\Outcross{{\beta'}_n}{{\alpha'}_{k+1}}{{x'}_{k+1}}{{\beta'}_k}{{\alpha'}_k}{{x'}_k}{{\alpha'}_1}{{x'}_1}\times\nonumber\\
&&\times\prod_{i\not=
k}\delta_{\beta_i,{\beta'}_i}\prod_{i\not=k,k+1}\delta_{\alpha_i,{\alpha'}_i}
\frac
{tr_q(\psi_{\beta_{k-1}\alpha_{k-1}}^{\beta_{k-2}}\psi_{\beta_{k}\alpha_{k}}^{\beta_{k-1}}\Rtff{\alpha_k}{\alpha_{k-1}}
\phi^{\beta_{k}\alpha_{k-1}}_{{\beta'}_{k-1}}\phi^{{\beta'}_{k-1}\alpha_{k}}_{\beta_{k-2}})}{[d_{\beta_{k-2}}]v_{{\beta}_k}^{{1\over
2}}v_{{\beta'}_k}^{-{1\over 2}}}
\delta_{\alpha_{k+1},{\alpha'}_k}\delta_{\alpha_{k},{\alpha'}_{k+1}}\\ &&{\bf
{\cal
A}^{elem}_{undercross.}}=\sum_{\beta_1,\cdots,\beta_n,{\beta'}_k}
\Incross{\beta_n}{\alpha_{k+1}}{x_{k+1}}{\beta_k}{\alpha_k}{x_k}{\alpha_1}{x_1}
\Outcross{{\beta'}_n}{{\alpha'}_{k+1}}{{x'}_{k+1}}{{\beta'}_k}{{\alpha'}_k}{{x'}_k}{{\alpha'}_1}{{x'}_1}\times\nonumber\\
&&\times\prod_{i\not=k}\delta_{\beta_i,{\beta'}_i}\prod_{i\not=k,k+1}\delta_{\alpha_i,{\alpha'}_i}
\frac
{tr_q(\psi_{\beta_{k-1}\alpha_{k-1}}^{\beta_{k-2}}\psi_{\beta_{k}\alpha_{k}}^{\beta_{k-1}}\Rtmff{\alpha_k}{\alpha_{k-1}}
\phi^{\beta_{k}\alpha_{k-1}}_{{\beta'}_{k-1}}\phi^{{\beta'}_{k-1}\alpha_{k}}_{\beta_{k-2}})}{[d_{\beta_{k-2}}]v_{{\beta}_k}^{{1\over
2}}v_{{\beta'}_k}^{-{1\over 2}}}
\delta_{\alpha_{k+1},{\alpha'}_k}\delta_{\alpha_{k},{\alpha'}_{k+1}}
\end{eqnarray} \ep

\medskip

\proof \\
The result for ${\cal A}^{elem}_{cup}$ and ${\cal A}^{elem}_{cap}$ is clearly
given by the Boltzmann
weight . The computation of the other elements need a careful description. The
idea is very simple.
We first absorb each link segment in the  attached boltzmann weight to put each
elements associated
to "square plaquettes" in a same practical form where all edges of the boundary
appear one and only
one time and always in the same order. This form allows us to reduce the gluing
of "plaquettes" to
one commutation plus one integration only. \\ In the case of an empty square
plaquette the
corresponding Boltzmann weight is already in the reduced form.\\
 For example the square plaquette element involved in the computation of a free
propagation is
given by:
\bea
&&\delta_{[{x'}_n,{x}_{n},{x}_{n-1},{x'}_{n-1}]}\guf{\alpha_{n-1}}_{[x_{n-1},{x'}_{n-1}]}=
\sum_{\beta_{n-1}{\beta'}_{n-2}}[d_{\beta_{n-1}}]\lambda_{\beta_{n-1}
\alpha_{n-1}\beta_{n-2}}^{-1}\times\\ &&\times
tr_{V_{\beta_{n-1}}}(\muf{\beta_{n-1}}\guf{\beta_{n-1}}_{[{x'}_{n},{x}_{n}]}
\guf{\beta_{n-1}}_{[x_{n},{x}_{n-1}]}\phi^{\beta_{n-1}
\alpha_{n-1}}_{{\beta'}_{n-2}}\guf{{\beta'}_{n-2}}_{[x_{n-1},{x'}_{n-1}]}\psi_{\beta_{n-1}
\alpha_{n-1}}^{{\beta'}_{n-2}}
\guf{\beta_{n-1}}_{[{x'}_{n-1},{x'}_{n}]})\nonumber
\eea where the
notations are summarized on the following figure:

\par
\centerline{\psfig{figure=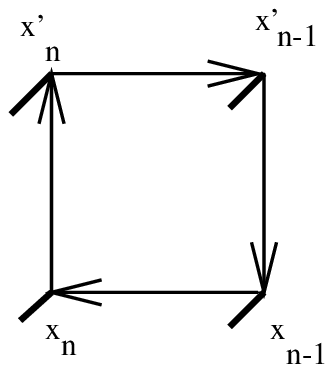}}
\par

In the case of a creation element for example, with a similar computation and
the same notations, we have :
\bea
&&\delta_{[{x'}_n,{x}_{n},{x}_{n-1},{x'}_{n-1}]}\guf{\alpha_{n-1}}_{[{x'}_{n},{x'}_{n-1}]}^{\#}=
\nonumber\\
&&=\sum_{\beta_n,{\beta'}_n} ([d_{\beta_n}][d_{{\beta'}_n}])^{1 \over 2}
N^{{\beta'}_n \alpha_n}_{\beta_n}tr_{V_{\beta_n}}(\muf{\beta_n}
\guf{\beta_n}_{[{x'}_n,x_n]}\guf{\beta_n}_{[x_n,x_{n-1}]}\guf{\beta_n}_{[x_{n-1},{x'}_{n-1}]}
\times\nonumber\\
&&\times\psi_{{\beta'}_n \alpha_n}^{\beta_n}
\guf{{\beta'}_n}_{[{x'}_n,{x'}_{n-1}]}\psi_{{\beta}_n
\bar{\alpha_n}}^{{\beta'}_n}\Rmff{\beta_n}{\bar{\alpha_n}}\lambda_{\beta_n
\bar{\alpha_n}
{\beta'}_n}) \eea

 When the square contains a crossing the reduction is less obvious and is given
by the
following lemma:
\bl[integration over a crossing]
With the notations of the figure:\\

\par
\centerline{\psfig{figure=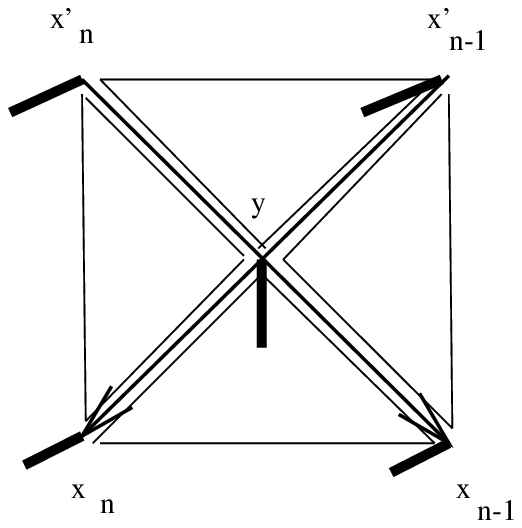}}
\par
it was shown in our last work \cite{BR2} that:
\bea
&&\int
dh(U_{[x_{n-1},y]})dh(U_{[y,{x'}_n]})dh(U_{[x_{n},y]})dh(U_{[y,{x'}_{n-1}]})\times\nonumber\\
&&\times\delta_{[y x_n x_{n-1}]}\delta_{[y x_{n-1} {x'}_{n-1}]}\delta_{[y
{x'}_{n-1}  {x'}_{n}]}
\delta_{[y {x'}_n x_{n}]}\guf{\alpha_n}_{[x_n y {x'}_{n-1}]}
\guf{\alpha_{n-1}}_{[x_{n-1} y
{x'}_n]}=\nonumber\\
&&=\sum_{\beta_n}[d_{\beta_n}](\frac{v_{{\beta'}_{n-1}}v_{{\beta}_{n-2}}}{v_{{\beta}_{n-1}}
v_{{\beta}_{n}}})^{1
\over 2}\frac
{tr_q(\psi_{\beta_{k-1}\alpha_{k-1}}^{\beta_{k-2}}\psi_{\beta_{k}\alpha_{k}}^{\beta_{k-1}}
\Rtff{\alpha_k}{\alpha_{k-1}}
\phi^{\beta_{k}\alpha_{k-1}}_{{\beta'}_{k-1}}
\phi^{{\beta'}_{k-1}\alpha_{k}}_{\beta_{k-2}})}{[d_{\beta_{k-2}}]}
tr_{V_{\beta_n}}(\muf{\beta_n} \guf{\beta_n}_{[{x'}_{n}
{x}_{n}]}\times\nonumber\\
&&\times\phi^{\beta_n \alpha_n}_{\beta_{n-1}}\guf{\beta_{n-1}}_{[{x}_{n}
{x}_{n-1}]}
\phi^{\beta_{n-1} \alpha_{n-1}}_{\beta_{n-2}}\guf{\beta_{n-2}}_{[{x}_{n-1}
{x'}_{n-1}]}
\psi^{\beta_{n-2}}_{{\beta'}_{n-1} \alpha_n}\guf{{\beta'}_{n-1}}_{[{x'}_{n-1}
{x'}_n]}\psi_{\alpha_{n-1}\beta_n
}^{{\beta'}_{n-1}}\Rpff{\alpha_{n-1}}{\beta_n}) \eea and the
analog relation for the undercrossing. \el

\medskip

Now all square plaquettes elements are in the reduced form.
Then the problem of computing the whole elementary block element is reduced to
the gluing of elements associated to each of the square plaquettes given in the
reduced form, i.e. a commutation $+$ an integration. Now, let us give a careful
computation
in the case of  ${\cal A}_{free}^{elem}$ and the other ones, very similar to
this one,
will be led to the reader.

The notations are summarized on the figure:

\par
\centerline{\psfig{figure=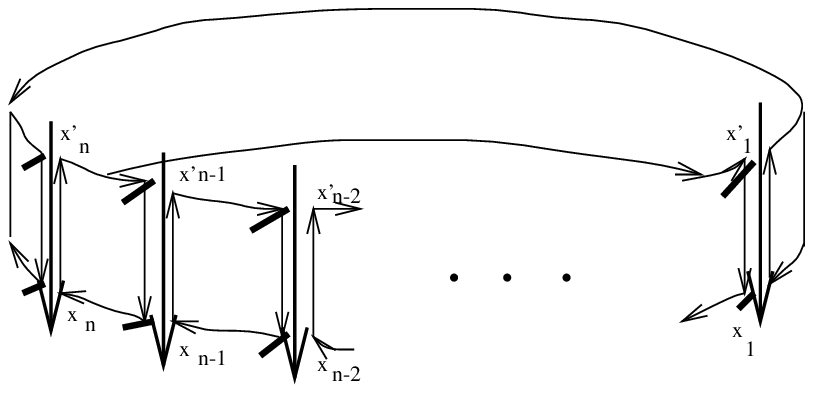}}
\par

({\bf Remark:} In the following computation we forget again the multiplicities
of
representations in all
decompositions, but we must take care of them...)
Let us first describe the gluing of two square plaquettes:

\ben
&&\int dh(U_{[x_{n-1},{x'}_{n-1}]})
\delta_{[{x'}_n,{x}_{n},{x}_{n-1},{x'}_{n-1}]}\guf{\alpha_{n-1}}_{[x_{n-1},{x'}_{n-1}]}
\delta_{[{x'}_{n-1},{x}_{n-1},{x}_{n-2},{x'}_{n-2}]}\guf{\alpha_{n-2}}_{[x_{n-2},{x'}_{n-2}]}=\\
&&=\int dh(U_{[x_{n-1},{x'}_{n-1}]})
\sum_{\beta_{n-1},\beta_{n-2},{\beta'}_{n-2},{\beta'}_{n-3}}[d_{\beta_{n-1}}][d_{\beta_{n-2}}]\lambda_{\beta_{n-1} \alpha_{n-1}\beta_{n-2}}^{-1}\lambda_{\beta_{n-2} \alpha_{n-2}\beta_{n-3}}^{-1}\times\\
&&\times
tr_{V_{\beta_{n-1}}}(\muf{\beta_{n-1}}\guf{\beta_{n-1}}_{[{x'}_{n},{x}_{n}]}\guf{\beta_{n-1}}_{[x_{n},{x}_{n-1}]}\phi^{\beta_{n-1} \alpha_{n-1}}_{{\beta'}_{n-2}}\guf{{\beta'}_{n-2}}_{[x_{n-1},{x'}_{n-1}]}\psi_{\beta_{n-1} \alpha_{n-1}}^{{\beta'}_{n-2}}
\guf{\beta_{n-1}}_{[{x'}_{n-1},{x'}_{n}]})\\
&&\times
tr_{V_{\beta_{n-2}}}(\muf{\beta_{n-2}}\guf{\beta_{n-2}}_{[{x'}_{n-1},{x}_{n-1}]}\guf{\beta_{n-2}}_{[x_{n-1},{x}_{n-2}]}\phi^{\beta_{n-2} \alpha_{n-2}}_{{\beta'}_{n-3}}\guf{{\beta'}_{n-3}}_{[x_{n-2},{x'}_{n-2}]}\psi_{\beta_{n-2} \alpha_{n-2}}^{{\beta'}_{n-3}}\guf{\beta_{n-2}}_{[{x'}_{n-2},{x'}_{n-1}]})=\\
&&=\sum_{(i),(j)}\int dh(U_{[x_{n-1},{x'}_{n-1}]})
\sum_{\beta_{n-1},\beta_{n-2},{\beta'}_{n-2},{\beta'}_{n-3}}[d_{\beta_{n-1}}][d_{\beta_{n-2}}]\lambda_{\beta_{n-1} \alpha_{n-1}\beta_{n-2}}^{-1}\lambda_{\beta_{n-2} \alpha_{n-2}{\beta'}_{n-3}}^{-1}\times\\
&&\times tr_{V_{\beta_{n-1}}\otimes
V_{\beta_{n-2}}}(\muf{\beta_{n-1}}\guf{\beta_{n-1}}_{[{x'}_{n},{x}_{n}]}\guf{\beta_{n-1}}_{[x_{n},{x}_{n-1}]}\phi^{\beta_{n-1} \alpha_{n-1}}_{{\beta'}_{n-2}}\muf{\beta_{n-2}}a_{(i)}^{\beta{n-2}}\muf{\beta_{n-2}}^{-1}\times\\
&&\guf{{\beta'}_{n-2}}_{[x_{n-1},{x'}_{n-1}]}
\muf{\beta_{n-2}}\guf{\beta_{n-2}}_{[{x'}_{n-1},{x}_{n-1}]}\guf{\beta_{n-2}}_{[x_{n-1},{x}_{n-2}]}
\phi^{\beta_{n-2}
\alpha_{n-2}}_{{\beta'}_{n-3}}\guf{{\beta'}_{n-3}}_{[x_{n-2},{x'}_{n-2}]}\psi_{\beta_{n-2} \alpha_{n-2}}^{{\beta'}_{n-3}}\times\\
&&\guf{\beta_{n-2}}_{[{x'}_{n-2},{x'}_{n-1}]}S(a_{(j)}^{\beta_{n-2}})
\psi_{\beta_{n-1} \alpha_{n-1}}^{{\beta'}_{n-2}}
b_{(i)}^{\beta_{n-1}}b_{(j)}^{\beta_{n-1}}\guf{\beta_{n-1}}_{[{x'}_{n-1},{x'}_{n}]})=\\
&&=\sum_{(i),(j)}\sum_{\beta_{n-1},\beta_{n-2},{\beta'}_{n-2},{\beta'}_{n-3}}[d_{\beta_{n-1}}][d_{\beta_{n-2}}]\lambda_{\beta_{n-1} \alpha_{n-1}\beta_{n-2}}^{-1}\lambda_{\beta_{n-2} \alpha_{n-2}{\beta'}_{n-3}}^{-1}\frac{\delta_{\beta_{n-2},{\beta'}_{n-2}}}{[d_{\beta_{n-2}}]}\times\\
&&\times
tr_{V_{\beta_{n-1}}}(\muf{\beta_{n-1}}\guf{\beta_{n-1}}_{[{x'}_{n},{x}_{n}]}\guf{\beta_{n-1}}_{[x_{n},{x}_{n-1}]}\phi^{\beta_{n-1} \alpha_{n-1}}_{{\beta}_{n-2}}\guf{\beta_{n-2}}_{[x_{n-1},{x}_{n-2}]}
\phi^{\beta_{n-2}
\alpha_{n-2}}_{{\beta'}_{n-3}}\guf{{\beta'}_{n-3}}_{[x_{n-2},{x'}_{n-2}]}\times\\
&&\times\psi_{\beta_{n-2}
\alpha_{n-2}}^{{\beta'}_{n-3}}\guf{\beta_{n-2}}_{[{x'}_{n-2},{x'}_{n-1}]}S(a_{(j)}^{\beta_{n-2}})
S^2(a_{(i)}^{\beta{n-2}})
\psi_{\beta_{n-1} \alpha_{n-1}}^{{\beta}_{n-2}}
b_{(i)}^{\beta_{n-1}}b_{(j)}^{\beta_{n-1}}\guf{\beta_{n-1}}_{[{x'}_{n-1},{x'}_{n}]})=\\
&&=\sum_{\beta_{n-1},\beta_{n-2},{\beta'}_{n-3}}[d_{\beta_{n-1}}]\lambda_{\beta_{n -1}\alpha_{n-1}\beta_{n-2}}^{-1}\lambda_{\beta_{n-2} \alpha_{n-2}{\beta'}_{n-3}}^{-1} tr_{V_{\beta_{n-1}}}(\muf{\beta_{n-1}}\guf{\beta_{n-1}}_{[{x'}_{n},{x}_{n}]}\guf{\beta_{n-1}}_{[x_{n},{x}_{n-1}]}\times\\
&&\phi^{\beta_{n-1}
\alpha_{n-1}}_{{\beta}_{n-2}}\guf{\beta_{n-2}}_{[x_{n-1},{x}_{n-2}]}
\phi^{\beta_{n-2}
\alpha_{n-2}}_{{\beta'}_{n-3}}\guf{{\beta'}_{n-3}}_{[x_{n-2},{x'}_{n-2}]}\psi_{\beta_{n-2} \alpha_{n-2}}^{{\beta'}_{n-3}}\guf{\beta_{n-2}}_{[{x'}_{n-2},{x'}_{n-1}]}\times\\
&&\times\psi_{\beta_{n-1}
\alpha_{n-1}}^{{\beta}_{n-2}}\guf{\beta_{n-1}}_{[{x'}_{n-1},{x'}_{n}]})\\
\een
In the same way we can glue $n-1$ Boltzmann weights and links. As a result we
have obviously:
\ben
&&\int \prod_{i=1}^{n-1} dh(U_{[x_{i},{x'}_{i}]})
\prod_{i=1}^{n-1}(\delta_{[{x'}_{i+1},{x}_{i+1},{x}_{i},{x'}_{i}]}\guf{\alpha_{i}}_{[x_{i},{x'}_{i}]})=\\
&&=\sum_{\beta_{n-1},\cdots,\beta_{2},{\beta}_{1}}[d_{\beta_{n-1}}]\prod_{i=2}^{n-1}\lambda_{\beta_i \alpha_{i}\beta_{i-1}}^{-1}\lambda_{\beta_{1} \alpha_{1}{\beta'}_{n}}^{-1} tr_{V_{\beta_{n-1}}}(\muf{\beta_{n-1}}\guf{\beta_{n-1}}_{[{x'}_{n},{x}_{n}]}\guf{\beta_{n-1}}_{[x_{n},{x}_{n-1}]}\phi^{\beta_{n-1} \alpha_{n-1}}_{{\beta}_{n-2}}\cdots\\
&&\cdots \phi^{\beta_{1}
\alpha_{1}}_{{\beta'}_{n}}\guf{{\beta'}_{n}}_{[x_{1},{x'}_{1}]}\psi_{{\beta}_{1} \alpha_{1}}^{{\beta'}_{n}}\cdots
\psi_{\beta_{n-1}
\alpha_{n-1}}^{{\beta}_{n-2}}\guf{\beta_{n-1}}_{[{x'}_{n-1},{x'}_{n}]})\\
\een
The computation of ${\cal A}^{elem}_{free}$ can be achieved by the gluing of
the last
Boltzmann weight to obtain the cylinder with $n$ strands.This Boltzmann weight
must be changed to the
square plaquette element corresponding to crossings,... in the computations of
the other blocks
elements. \ben &&\int \prod_{i=1}^{n} dh(U_{[x_{i},{x'}_{i}]})
\prod_{i=1}^{n}(\delta_{[{x'}_{i+1},{x}_{i+1},{x}_{i},{x'}_{i}]}\guf{\alpha_{i}}_{[x_{i},{x'}_{i}]})=\\
&&=\int dh(U_{[x_{n},{x'}_{n}]})dh(U_{[x_{1},{x'}_{1}]})\times\\
&&(\sum_{\beta_n,{\beta'}_{n-1}}[d_{\beta_n}]\lambda_{\beta_n \alpha_n
{\beta'}_{n-1}}^{-1}tr_{V_{\beta_n}}(\muf{\beta_n}\guf{\beta_n}_{[x_1,x_n]}
\phi^{\beta_n
\alpha_n}_{{\beta'}_{n-1}}\guf{{\beta'}_{n-1}}_{[x_n,{x'}_n]}\psi_{\beta_n
\alpha_n}^{{\beta'}_{n-1}}\guf{\beta_n}_{[{x'}_n,{x'}_1]}\guf{\beta_n}_{[{x'}_1,{x}_1]})) \times\\
&&\times
(\sum_{\beta_{n-1},\cdots,\beta_{2},{\beta}_{1}}[d_{\beta_{n-1}}]\prod_{i=2}^{n-1}\lambda_{\beta_i
\alpha_{i}\beta_{i-1}}^{-1}\lambda_{\beta_{1} \alpha_{1}{\beta'}_{n}}^{-1}
tr_{V_{\beta_{n-1}}}(\muf{\beta_{n-1}}\guf{\beta_{n-1}}_{[{x'}_{n},{x}_{n}]}\guf{\beta_{n-1}}_{[x_{n},{x}_{n-1}]}
\phi^{\beta_{n-1} \alpha_{n-1}}_{{\beta}_{n-2}}\cdots\\ &&\cdots
\phi^{\beta_{1}
\alpha_{1}}_{{\beta'}_{n}}\guf{{\beta'}_{n}}_{[x_{1},{x'}_{1}]}\psi_{{\beta}_{1}
\alpha_{1}}^{{\beta'}_{n}}\cdots \psi_{\beta_{n-1}
\alpha_{n-1}}^{{\beta}_{n-2}}\guf{\beta_{n-1}}_{[{x'}_{n-1},{x'}_{n}]}))=\\
&&=\int
dh(U_{[x_{n},{x'}_{n}]})dh(U_{[x_{1},{x'}_{1}]})
\sum_{{\beta'}_{n-1},\beta_n,\cdots,\beta_{2},{\beta'}_{1}}[d_{\beta_{n-1}}][d_{\beta_n}]
\prod_{i=2}^{n-1}\lambda_{\beta_i \alpha_{i}\beta_{i-1}}^{-1}\lambda_{\beta_{1}
\alpha_{1}{\beta'}_{n}}^{-1}  \lambda_{\beta_n \alpha_n
{\beta'}_{n-1}}^{-1}\times\\
&&tr_{V_{\beta_n} \otimes
V_{\beta_{n-1}}}(\muf{\beta_n}\guf{\beta_n}_{[x_1,x_n]} \phi^{\beta_n
\alpha_n}_{{\beta'}_{n-1}} \muf{\beta_{n-1}} a_{(j)}^{\beta_{n-1}}
\muf{\beta_{n-1}}^{-1}
\guf{{\beta'}_{n-1}}_{[x_n,{x'}_n]}\muf{\beta_{n-1}}\guf{\beta_{n-1}}_{[{x'}_{n},{x}_{n}]}
\guf{\beta_{n-1}}_{[x_{n},{x}_{n-1}]}\times\\  &&\times\phi^{\beta_{n-1}
\alpha_{n-1}}_{{\beta}_{n-2}}\cdots\guf{\beta_{1}}_{[x_{1},{x}_{n}]}b_{(i)}^{\beta_1}\phi^{\beta_{1}
\alpha_{1}}_{{\beta'}_{n}} \psi_{\beta_n
\alpha_n}^{{\beta'}_{n-1}}b_{(j)}^{\beta_n}\guf{{\beta}_{n}}_{[{x'}_n,{x'}_1]}
\guf{\beta_n}_{[{x'}_1,{x}_1]}\guf{{\beta'}_{n}}_{[x_{1},{x'}_{1}]}S(a_{(i)}^{{\beta'}_{n}})
\psi_{{\beta}_{1} \alpha_{1}}^{{\beta'}_{n}}\cdots\\ &&\cdots\psi_{\beta_{n-1}
\alpha_{n-1}}^{{\beta}_{n-2}}\guf{\beta_{n-1}}_{[{x'}_{n-1},{x'}_{n}]}))=\\
&&=\sum_{\beta_n,\cdots,{\beta}_{1}}[d_{\beta_n}]
\prod_{i=1}^{n}v_{\alpha_{i}}^{-{1 \over
2}}tr_{V_{\beta_n}}(\muf{\beta_n}\guf{\beta_n}_{[x_1,x_n]} \phi^{\beta_n
\alpha_n}_{{\beta}_{n-1}}\guf{\beta_{n-1}}_{[x_{n},{x}_{n-1}]}\phi^{\beta_{n-1}
\alpha_{n-1}}_{{\beta}_{n-2}}\cdots
\guf{\beta_{1}}_{[x_{1},{x}_{n}]}b_{(i)}^{\beta_1}\phi^{\beta_{1}
\alpha_{1}}_{{\beta}_{n}} \times\\
&&\times\muf{{\beta}_{n}}^{-1}S(a_{(i)}^{{\beta}_{n}}))
tr_{V_{{\beta}_{n-1}}}(S^2(a_{(j)}^{{\beta}_{n-1}})\psi_{\beta_n
\alpha_n}^{{\beta}_{n-1}}
b_{(j)}^{\beta_n} \guf{{\beta}_{n}}_{[{x'}_n,{x'}_1]} \psi_{{\beta}_{1}
\alpha_{1}}^{{\beta'}_{n}}\cdots \psi_{\beta_{n-1}
\alpha_{n-1}}^{{\beta}_{n-2}}\guf{\beta_{n-1}}_{[{x'}_{n-1},{x'}_{n}]})=\\
&&=\sum_{\beta_1,m_1,\cdots,\beta_n,m_n} \Intri{\beta_n
m_n}{\alpha_{n}}{x_{n}}{\beta_1
m_1}{\alpha_1}{x_1} \Outtri{{\beta}_n m_n}{{\alpha}_{n}}{{x}_{n}}{{\beta}_1
m_1}{{\alpha}_1}{{x'}_1}\\ \een This concludes the computation of ${\cal
A}^{elem}_{overcross}$,
${\cal A}^{elem}_{undercross}$, ${\cal A}^{elem}_{creation}$, ${\cal
A}^{elem}_{annihil}$, and
${\cal A}^{elem}_{free}.$\\ The computations of ${\cal A}^{elem}_{(n,m)(n+m)
tri.}$ and ${\cal
A}^{elem}_{(n,m)(n+m) tri.}$ need one more step. It uses naturally the
expression of ${\cal
A}^{elem}_{free}$ as a basic object. Indeed we compute the element associated
to the trinion by
gluing one more plaquette to the cylinder with $n$ strands as it is shown in
the following figure.

\par
\centerline{\psfig{figure=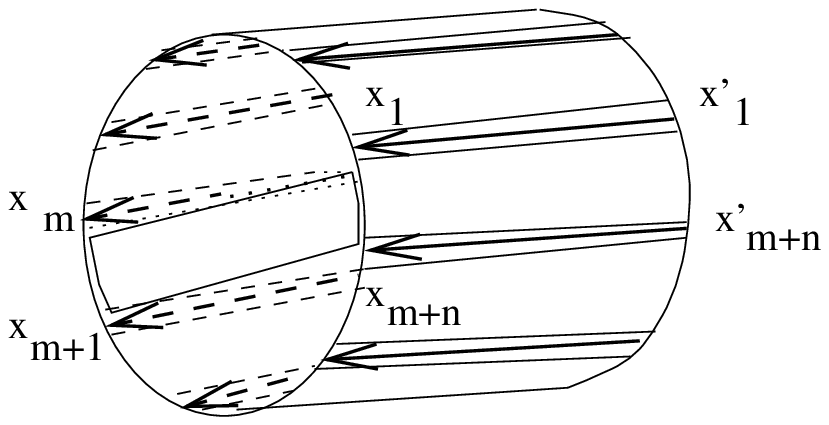}}
\par

\medskip
the computation is realized by the usual techniques:\\
\ben
&&{\cal A}_{(n,m)(n+m) tri.}=\int dh(U_{[x_1,x_{m+n}]})dh(U_{[x_{m+1},x_{m}]})
\delta_{[x_1,x_m,x_{m+1},x_{m+n}]}\times\\
&&\times\sum_{\beta_1,\cdots,\beta_{m+n}}
\Intri{\beta_{m+n}}{\alpha_{n+m}}{x_{n+m}}{\beta_1}{\alpha_1}{x_1}\Outtri{{\beta}_{m+n}}{{\alpha}_{n+m}}{{x'}_{n+m}}{{\beta}_1}{{\alpha}_1}{{x'}_1}=\\
&&\int dh(U_{[x_1,x_{m+n}]})dh(U_{[x_{m+1},x_{m}]})
\sum_{{\beta'}_{m+n}}[d_{{\beta'}_{m+n}}]v_{{\beta'}_{m+n}}^{-1}
tr_{V_{{\beta'}_{m+n}}}(\muf{{\beta'}_{m+n}}\guf{{\beta'}_{m+n}}_{[x_{m+n},x_1]}
\guf{{\beta'}_{m+n}}_{[x_{1},x_m]}\times\\
&&\times\guf{{\beta'}_{m+n}}_{[x_{m},x_{m+1}]}
\guf{{\beta'}_{m+n}}_{[x_{m+1},x_{m+n}]})\;
\sum_{\beta_{m+n},\cdots,{\beta}_{1}}
tr_{V_{\beta_{n+m}}}(\muf{\beta_{n+m}}\Rff{\beta_{m+n}}{\alpha_1}
\guf{\beta_{n+m}}_{[x_1,x_{n+m}]}\phi^{\beta_{m+n}\alpha_{m+n}}_{{\beta}_{m+n-1}}\cdots\\
&&\cdots\guf{\beta_{1}}_{[x_{2},{x}_{1}]}\phi^{\beta_{1}
\alpha_{1}}_{{\beta}_{m+n}})
\;\Outtri{{\beta}_{m+n}}{{\alpha}_{n+m}}{{x'}_{n+m}}{{\beta}_1}{{\alpha}_1}{{x'}_1}=\\
&&=\int dh(U_{[x_1,x_{m+n}]})dh(U_{[x_{m+1},x_{m}]})
\sum_{{\beta'}_{m+n},\beta_{m+n},\cdots,{\beta}_{1}}
[d_{{\beta'}_{m+n}}]v_{{\beta'}_{m+n}}^{-1}\times\\
&&\times tr_{V_{{\beta'}_{m+n}}\otimes V_{\beta_{n+m}}}
(\muf{{\beta'}_{m+n}}\muf{\beta_{n+m}}\Rff{\beta_{n+m}}{\alpha_1}
S^{-1}(b^{\beta_{m+n}}_{(j)})\muf{\beta_{n+m}}^{-1}
(\guf{{\beta'}_{m+n}}_{[x_{m+n},x_1]})\muf{\beta_{n+m}}\guf{\beta_{n+m}}_{[x_1,x_{n+m}]})\times\\
&&\times a_{(j)}^{{\beta'}_{n+m}}a_{(i)}^{\beta_{n+m}}
\guf{{\beta'}_{m+n}}_{[x_{1},x_m]}
\phi^{\beta_{m+n}\alpha_{m+1}}_{{\beta}_{m+n-1}}a_{(l)}^{{\beta}_{n+m-1}}\cdots
\phi_{\beta_{m}}^{\beta_{m+1}\alpha_{m+1}}
S^{-1}(b^{\beta_{m}}_{(k)})\muf{\beta_{m+1}}^{-1}\times\\
&&\times(\guf{{\beta'}_{m+n}}_{[x_{m},x_{m+1}]}\muf{\beta_{m+1}}\guf{{\beta}_{m+1}}_{[x_{m+1},x_{m}]})
a_{(k)}^{{\beta'}_{n+m}}\guf{{\beta'}_{m+n}}_{[x_{m+1},x_{m+n}]}b_{(i)}^{{\beta'}_{m+n}}
S(b^{{\beta'}_{m+n}}_{(l)})
\phi_{\beta_{m-1}}^{\beta_m \alpha_m}\cdots\\
&&\cdots\guf{\beta_{1}}_{[x_{2},{x}_{1}]}\phi^{\beta_{1}
\alpha_{1}}_{{\beta}_{m+n}})
\Outtri{{\beta}_{m+n}}{{\alpha}_{n+m}}{{x'}_{n+m}}{{\beta}_1}{{\alpha}_1}{{x'}_1}=\\
&&=\sum_{{\beta'}_{m+n},\beta_{m+n},\cdots,{\beta}_{1}}v_{\beta_m}
[d_{{\beta}_{m+n}}]^{-1} \delta_{\beta_{m+n},\beta_m,{\beta'}_{m+n}}\times\\
&&\times tr_{V_{{\beta}_{m}}}
(\muf{{\beta}_{m}}\muf{\beta_{n+m}}\Rff{\beta_{m}}{\alpha_1}
\guf{{\beta}_{m}}_{[x_{1},x_m]}
\phi^{\beta_{m}\alpha_{m}}_{{\beta}_{m-1}}\cdots\guf{{\beta}_{1}}_{[x_{2},x_1]}
\phi_{\beta_{m}}^{\beta_{1}\alpha_{1}})\times\\
&&\times tr_{V_{{\beta}_{m+n}}}(\muf{\beta_{m+n}}a_{(i)}^{{\beta}_{n+m}}
\phi_{\beta_{m+n-1}}^{\beta_{m+n} \alpha_{m+n}}a_{(l)}^{{\beta}_{n+m-1}}\cdots
\phi^{\beta_{m+1}
\alpha_{m+1}}_{{\beta}_{m+n}}\guf{{\beta}_{m+n}}_{[x_{m+1},x_{m+n}]}
S(b^{{\beta}_{m+n}}_{(l)})b_{(i)}^{{\beta}_{m+n}})\times\\
&&\times\Outtri{{\beta}_{m+n}}{{\alpha}_{n+m}}{{x'}_{n+m}}{{\beta}_1}{{\alpha}_1}{{x'}_1}=\\
&&= \sum_{\beta_1,\cdots, \beta_{n+m}} [d_{\beta_{n+m}}]^{-1}
\Intri{{\beta'}_n}{{\alpha'}_n}{{x'}_n}{{\beta'}_1}{{\alpha'}_1}{{x'}_1}\Intri{{\beta''}_n}{{\alpha''}_n}{{x''}_n}{{\beta''}_1}{{\alpha''}_1}{{x''}_1}
\times\\ &&\times
\Outtri{{\beta}_n}{{\alpha}_n}{{x}_n}{{\beta}_1}{{\alpha}_1}{{x}_1}
\delta_{{\beta'}_n, {\beta''}_m, {\beta}_{n+m}, \beta_m}
\prod_{k=1}^{n}\delta_{\alpha_{m+k}, {\alpha'}_{k}} \prod_{k=1}^{m}
\delta_{\alpha_{k}, {\alpha''}_{k}}
\prod_{k=1}^{n-1}\delta_{\alpha_{m+k}, {\alpha'}_{k}} \prod_{k=1}^{m-1}
\delta_{\alpha_{k}, {\alpha''}_{k}}\nonumber\\
&&=\sum_{\beta_1,\cdots, \beta_{n+m}} [d_{\beta_{n+m}}]^{-1}
\Intri{{\beta}_n}{{\alpha}_n}{{x}_n}{{\beta}_1}{{\alpha}_1}{{x}_1}
\Outtri{{\beta'}_n}{{\alpha'}_n}{{x'}_n}{{\beta'}_1}{{\alpha'}_1}{{x'}_1}\times\\
&&\times
\Outtri{{\beta''}_n}{{\alpha''}_n}{{x''}_n}{{\beta''}_1}{{\alpha''}_1}{{x''}_1}
\delta_{{\beta'}_n, {\beta''}_m, {\beta}_{n+m}, \beta_m}
\prod_{k=1}^{n}\delta_{\alpha_{m+k}, {\alpha'}_{k}} \prod_{k=1}^{m}
\delta_{\alpha_{k}, {\alpha''}_{k}}
\prod_{k=1}^{n-1}\delta_{\alpha_{m+k}, {\alpha'}_{k}}
\prod_{k=1}^{m-1}\delta_{\alpha_{k}, {\alpha''}_{k}}\nonumber\\
\een
This ends the proof of the Theorem.\\

\cqfd

These results can be easily generalized to the case where $q$ is a root of
unity.
 It suffices to realize the following replacements:\\
The expression for the In state becomes
\begin{eqnarray*}
&&\Instate{\beta_n}{\alpha_n}{x_n}{\beta_{n-1}}{\alpha_{n-1}}{x_{n-1}}{\beta_1}{\alpha_1}{x_1}=\\
&&=v_{\beta_n}tr_{V_{\beta_n}}(S(A)\muf{\beta_n}
{\buildrel {\alpha_1} \over {\theta^{(1)}_l}}
{\buildrel {\beta_n} \over {\theta^{(2)}_l}}
\Rff{\beta_n}{\alpha_1}
{\buildrel {\beta_n} \over {\theta^{-1 (1)}_k}}
{\buildrel {\alpha_1} \over {\theta^{-1 (2)}_k}}
\guf{\beta_n}_{[x_1,x_{n}]}
S({\buildrel {\beta_n} \over {\theta^{(1)}_m}}){\buildrel {\beta_n} \over {A}}
{\buildrel {\beta_n} \over {\theta^{(2)}_m}}\otimes {\buildrel {\alpha_n} \over
{\theta^{(3)}_m}}
\phi^{\beta_n \alpha_{n}}_{\beta_{n-1}} \cdots\\  &&\cdots
\phi^{\beta_2 \alpha_2}_{\beta_1}\guf{\beta_1}_{[x_2,x_{1}]} S({\buildrel
{\beta_1} \over
{\theta^{(1)}_i}}) S({\buildrel {\beta_1} \over {\theta^{(1)}_j}}){\buildrel
{\beta_1} \over
{A}}{\buildrel {\beta_1} \over {\theta^{(2)}_j}}\otimes {\buildrel {\alpha_1}
\over
{\theta^{(3)}_j}} \phi^{\beta_1 \alpha_1}_{\beta_n} {\buildrel {\beta_n} \over
{\theta^{(2)}_i}}
{\buildrel {\beta_n} \over {B}}
S({\buildrel {\beta_n} \over {\theta^{(3)}_i}})
S({\buildrel {\beta_n} \over {\theta^{-1 (3)}_k}})
S({\buildrel {\beta_n} \over {\theta^{(3)}_l}}))
\end{eqnarray*}
and the analog formula for the Out state.
The $6-j$ associated to the crossing is changed to obtain an intertwiner.\\
The other objects remain unchanged and the summations are restricted to
physical representations
only.\\

\medskip

We have the following corollary:
\bp
If $L$ is a link without boundaries in $D\times[0,1]$ we have for any value of
$q$:
\be
\frac{<W_L>_{q-YM(S^2)}}{<1>_{q-YM(S^2)}}= RT_{{\cal U}_q({\cal G})}(L)
\ee
where $RT$ is the Reshetikhin-Turaev's quantum invariant of coloured links.
This result has already been shown in \cite{BR2}.
and generally if $\Sigma$ is a closed surface, the invariant associated to
$\Sigma+L$ is simply
a generalization to the case of a surface of the Reshetikhin-Kirillov invariant
in the shadow world \cite{KR}.
This theorem can be considered as a proof of the equivalence of Invariants
arising from Chern-Simons theory and Reshetikhin-Turaev quantum invariants.
\ep

\bigskip

\bd
We can also generalize our invariant to admit other objects called "coupons"
defined to be respectively represented on the following figure:

\par
\centerline{\psfig{figure=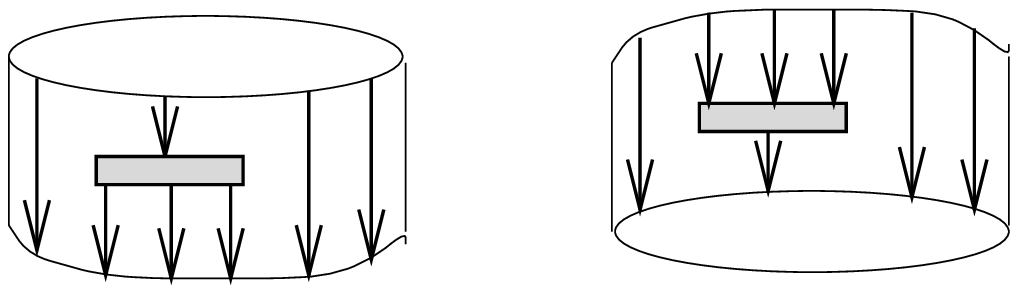}}
\par

and which expressions are given by:
\ben
{\cal
A}_{(n+k)(n)coupon}=\sum_{}\Intri{\beta_{n+k}}{\alpha_{n+k}}{x_{n+k}}{\beta_1}{\alpha_1}{x_1}\times
Shadow(coupon)\times
\Outtri{{\beta'}_{n}}{{\alpha'}_{n}}{{x'}_{n}}{\beta_1}{\alpha_1}{x_1}\\
{\cal
A}_{(n)(n+k)coupon}=\sum_{}\Intri{\beta_{n}}{\alpha_{n}}{x_{n}}{\beta_1}{\alpha_1}{x_1}\times
Shadow(coupon)\times
\Outtri{{\beta'}_{n+k}}{{\alpha'}_{n+k}}{{x'}_{n+k}}{\beta_1}{\alpha_1}{x_1}
\een
with $Shadow(coupon)$ being defined as usual by the following rule:
\par
\centerline{\psfig{figure=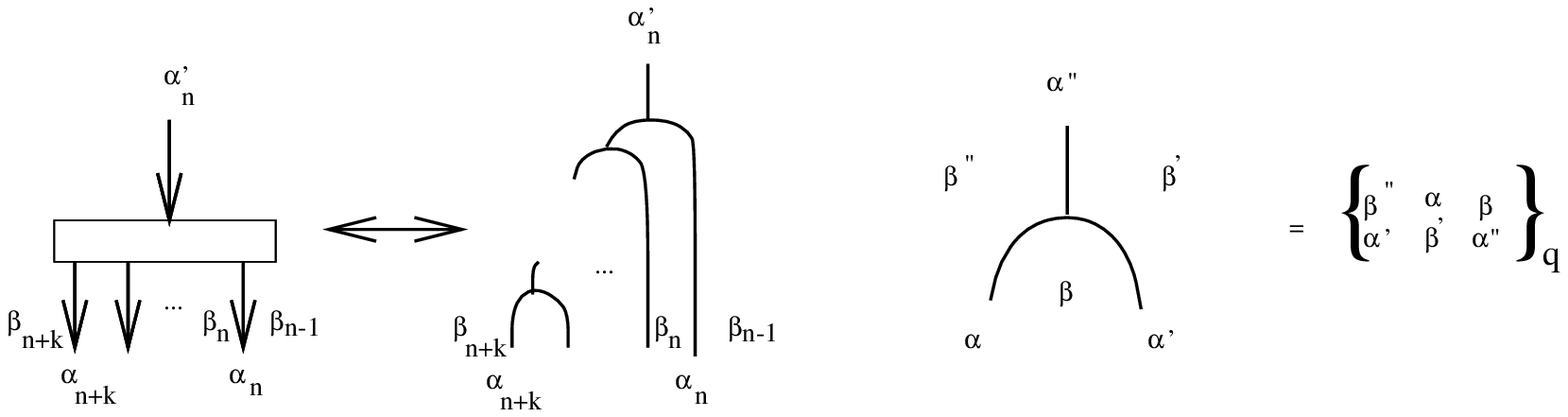}}
\par
\ed

\bigskip

\section{A new description of invariants of three manifolds}

In this chapter $q$ will be a root of unity.

\subsection{Heegaard splitting and surgery of 3-manifolds}
 In the following ${\cal M}$ is a compact orientable 3-manifold given
by a simplicial complex $K$.Let us recall standard definitions that can be
found in \cite{Si}.
\bd
A {\bf canonical region} ${\cal R}$ of ${\cal M}$ is a region within which
there are p non intersecting 2-cells $(E_i)_{i=1 \cdots p}$(the {\bf canonical
cells}) with boundaries $e_i$ (the {\bf canonical curves}) on the boundary
${\cal L}$ of ${\cal R}$ such that we obtain a 3-cell by cutting ${\cal R}$ at
each $(E_i).$
 A surface ${\cal L}$ is said to be a {\bf canonical surface} of a 3-manifold
${\cal M}$ if it satisfies these conditions:
\begin{itemize}
\item ${\cal L}$ is a subcomplex of  ${\cal M}$ and is a compact, connected
2-dimensional manifold
\item  ${\cal M}= {\cal R}_1 + {\cal L} + {\cal R}_2$ with ${\cal R}_1,{\cal
R}_2$ canonical regions and
${\cal L}=\partial {\cal R}_1 = \partial {\cal R}_2$
\end{itemize}
such a decomposition is called {\bf canonical decomposition}.
It is important to recall that, if $g$ is the genus of ${\cal L}$, in this
case, ${\cal R}_1$ and ${\cal R}_2$ are homeomorphic to a genus $g$ handlebody.
\ed

{\bf Remark: }
It is easy to give, for each  ${\cal M}$, at least one  canonical
decomposition.
To this aim, let us consider $\{A_0^i\}$, $\{A_1^j\}$, $\{A_2^k\}$, $\{A_3^l\}$
the sets of $0-,1-,2-,3-$cells of $K$.
Let $\{B_1^j\}$, $\{B_2^k\}$, $\{B_3^l\}$ be respectively the middle of
$\{A_1^j\}$, $\{A_2^k\}$, $\{A_3^l\}$. The complex $K'$ obtained by adding
the $B_{1}s$, the $B_{2}s$ and the $B_{3}s$ to the vertices of $K$ is called
the {\bf first derived complex of K}, its 3-simplexes are of the form
$(A_{0}B_{1}B_{2}B_{3})$.
The {\bf second derived complex of K}, denoted $K"$ is the complex complex
generated from $K'$ by adding the vertices $C_{1}s,C_{2}s,C_{3}s$ middle of
the $1-,2-,3-$simplexes of $K'$.
Let us denote by ${\cal R}_1$ the set of all $3-$simplexes of $K''$ of the type
$(A_{0}C_{1}C_{2}C_{3})$ or $(B_{1}C_{1}C_{2}C_{3})$, by ${\cal R}_2$ the set
of
all $3-$simplexes of $K"$ of the type $(B_{2}C_{1}C_{2}C_{3})$ or
$(B_{3}C_{1}C_{2}C_{3})$ and by ${\cal L}$ the common frontier of ${\cal R}_1$
and ${\cal R}_2$.
If we call $G$ (resp.$G^{\star}$) the linear graph
generated by the 1-simplexes of $K$ (resp. its dual) we can see that ${\cal
R}_1$
and ${\cal R}_2$ are respectively $K''-$neighbourhood of $G$ and $G^{\star}$.
Then we have that ${\cal M}= {\cal R}_1 + {\cal L} + {\cal R}_2$ is a canonical
decomposition, it is
called the {\bf canonical decomposition derived from the triangulation}.
It is easy to check that for each canonical decomposition, there exits a
triangulation of ${\cal M}$ such that the decomposition is in fact the
canonical decomposition derived from the triangulation.
\bd {\bf Heegaard Splitting}\\
A Heegaard splitting of a 3-manifold ${\cal M}$ is a set $(g,f)$ where $g$
is a non negative integer and $f$ is a diffeomorphism of a genus $g$ surface
${\cal L}_g$ such that ${\cal M}$ is the manifold obtained by gluing two copies
of the handlebody ${\cal T}_g$ (the interior of ${\cal L}_g$) along their
boundaries after
having acted on one of them by f:
$${\cal M} = {\cal T}_g \#_f {\cal T}_g$$
A Heegaard diagram  is a set $({\cal L},(e_i)_{i=1\cdots g},(f_j)_{j=1 \cdots
g})$ where ${\cal L}$ is a compact connected 2-dimensional manifold of genus
$g$ and $(e_i)_{i=1\cdots g}$ (resp.$(f_j)_{j=1 \cdots g})$) are canonical
curves of ${\cal R}_1$, the region interior to  ${\cal L}$ (resp. canonical
curves of ${\cal R}_2$ the
exterior of ${\cal L}$).This data is sufficient to reconstruct an element $f$
of $Diff({\cal L})$ such that $f(e_j)=f_j.$
Two Heegaard diagrams are said to be equivalent if they describe homeomorphic
3-manifolds.
Let $({f'}_i)_{i=1...g}$ be $g$ other canonical curves in ${\cal L}$ such that
$({\cal L},(e_i)_{i=1\cdots g},({f'}_j)_{j=1 \cdots g})$ is a Heegaard diagram
of the sphere $S^3$ then $({\cal L},(e_i)_{i=1\cdots g},(f_j)_{j=1 \cdots
g},({f'}_j)_{j=1 \cdots g})$ is said to be an augmented Heegaard diagram.
\ed
We must recall that any element of the moduli space of a surface can be written
as the composition of Dehn twists. A Dehn twist can be described by the
following replacement of a regular neighbourhood of the corresponding curve:

\par
\centerline{\psfig{figure=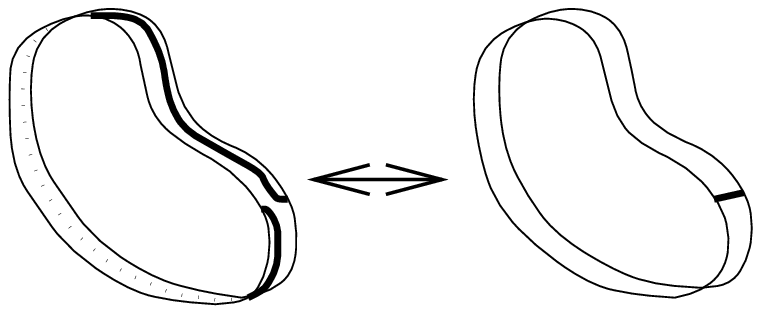}}
\par

There is an important theorem due to Singer \cite{Si} describing the
relation between equivalent Heegaard diagrams.

\bd{\bf Singer's elementary moves}
Let us describe a set of elementary moves on the Heegaard diagrams:
\begin{itemize}
{\bf type 0: trivial moves}
\item replace a curve by another curve isotopic to it, or to its inverse, or
reembedded the canonical surface ${\cal L}$ in a different way in $S^3$.\\

{\bf type 1: solid handlebodies diffeomorphisms}
\item replace one canonical curve of the set $(e_i)_{i=1\cdots g}$ (resp.
$(f_j)_{j=1\cdots g}$ ) by the composition
of this curve with another one in this set.
\item making a Dehn twist along one of the $e_{i}s$.\\

{\bf type 2: $g \rightarrow g+1$ moves}
\item add a handle to ${\cal L}$ and define $e_{p+1}$ (resp. $f_{p+1}$) to be
the $a-$cycle (resp. the $b-$cycle) of this handle, or erase a handle with
cycles for which $e_{p+1}$ is the $a-$cycle (resp. $f_{p+1}$ is the $b-$cycle).
\end{itemize}
\ed

Then we have the following classification theorem \cite{Si}:

\bp[ Singer's Theorem]
If the diagrams $D$ and $D'$, related by a finite number of Singer's moves,
give rise to the manifolds ${\cal M}$ and ${\cal M}'$ then ${\cal M}$ and
${\cal M}'$ are homeomorphic.
Conversely, if $D$ and $D'$ are any two Heegaard diagrams whatsoever arising
from a manifold ${\cal M}$ then $D$ and $D'$ are related by a finite number of
Singer's moves.
\ep
A more generally used description of three manifolds is "the surgery
presentation". Let us recall some facts about this description \cite{Li}.
\bd[Surgery presentation of 3-manifolds]
Let $(R,r)=\cup_{i=1}^{n}(R_i,r_i)$ be a framed link in the oriented sphere
$S^3.$ We can define a manifold ${\cal M}$
by "surgery" from $(R,r)$ using the following procedure:\\
remove from $S^3$ pairwise disjoint tubular neighbourhoods $V_i$ of the curves
$R_i$ and resew them identifying a meridian $z_i$ in $\partial V_i$ with a
curve $y_i \in \partial(S^3\setminus V_i^{int})$  which links $R_i$ exactly
$r_i$ times.\\
Moreover, every 3-manifold ${\cal M}$ can be obtained from a certain framed
link  by this procedure \cite{Li}.
\ed
It is relatively easy to relate the Heegaard and Surgery points of view
\cite{NL}. Let us consider a Heegaard diagram based on a gluing diffeomorphism
$f$ described in terms of Dehn twists of the surface. We first remark that
splitting $S^3$ along ${\cal L}$ then doing a Dehn twist along a certain curve
and resewing the handlebody is equivalent to do a surgery along the ribbon
glued on the surface along this curve as it can be seen on the figure:

\par
\centerline{\psfig{figure=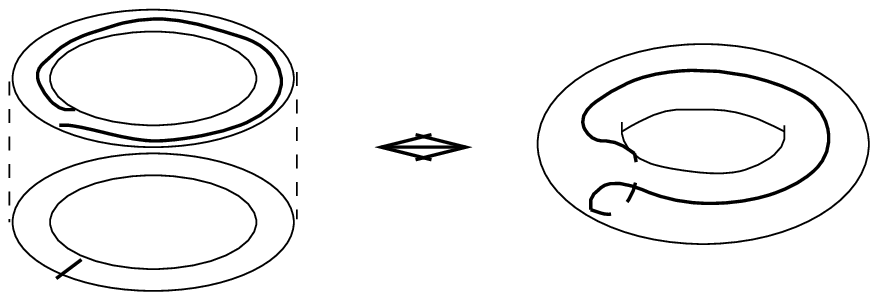}}
\par

Let $(R_i)_{i=1..n}$ be a set of ribbons trivially embedded on the surface
${\cal L}$, $f_i$ the corresponding Dehn twists. We want also define the framed
link $L$ defined to be the set of ribbons $(R_i \times \epsilon_i)_{i=1 \cdots
n}$ for $0 \le \epsilon_1 \le \cdots \le \epsilon_n \le 1.$
We consider a partition of $S^3$ in three pieces: ${\cal L}\times [0,1]$, the
handlebody ${\cal H}_g$ interior to ${\cal L}\times \{0\}$ and the handlebody
${\cal H}^{'}_g$ exterior to ${\cal L}\times \{1\}.$ Let us consider the
manifold ${\cal M}(f_1,f_2,\cdots,f_n)$ obtained by gluing the manifolds ${\cal
H}_g$, ${\cal H}^{'}_g$, $({\cal L}\times [\epsilon_{i-1},\epsilon_i])_{i=1
\cdots n}$ with the gluing diffeomorphisms $id, f_1, \cdots, f_n.$ Obviously
the manifold  ${\cal M}(f_1,f_2,\cdots,f_n)$ is the manifold defined by the
Heegaard data $({\cal L}, f_n \circ f_{n-1} \circ \cdots \circ f_1)$, but it is
also obvious that this manifold is that defined by the surgery data $R.$ We
will say that these surgery and Heegaard presentation of the same manifold are
"related" description of ${\cal M}.$

\subsection{Invariants associated to Heegaard diagrams and Lattice q-gauge
theory}

Our principal aim in this section is to prove the following theorem:

\bp[Invariants of three manifolds and Heegaard diagrams]
Let $({\cal T}_g, (x_j)_{j=1,...,g}, (y_j)_{j=1,...,g}, (z_j)_{j=1,...,g})$ be
an augmented Heegaard diagram associated to a manifold ${\cal M}$ then the
expectation value :
\be
{\cal J}_{{\cal M}}=\frac{<
\prod_{i=1}^{g}\delta_{y_i}\prod_{i=1}^{g}\delta_{x_i}>_{q-YM({\cal T}_g)}}{<
\prod_{i=1}^{g}\delta_{z_i}\prod_{i=1}^{g}\delta_{x_i}>_{q-YM({\cal T}_g)}}
\ee
 is an invariant of the manifold ${\cal M}.$ Moreover this value is equal to
the Reshetikhin-Turaev invariant associated to the manifold ${\cal M}.$
\ep
The normalization by the expectation value associated to the sphere is chosen
to obtain an $3-$manifold invariant equal to $1$ for the sphere.
{\bf Remark 1:}
The latter definition of the correlation function  is in fact very natural from
the general construction of q-gauge theory. Indeed putting a delta function
associated to a plaquette $P$ corresponds to imposing that any ribbon, i.e.
holonomy defined in terms of the gauge fields algebra, can  be displaced
through $P$ without torsion and without changing the expectation value. So,
adding to the projector associated to the surface some
delta functions corresponding to the $x_i$s, i.e. canonical curves of the
interior handlebody, and, at a future time, the delta functions of the $y_i$s,
i.e. canonical curves of the exterior handlebody, allows us to displace any
curve through a handle of any of the two Heegaard components.This is exactly
what we want to do in the framework of Chern-Simons theory.\\
{\bf Remark 2:}
Using the properties
\bea
&&(\frac{\delta_{C}}{\sum_{\alpha}[d_{\alpha}]^2})^2=(\frac{\delta_{C}}{\sum_{\alpha}[d_{\alpha}]^2})\\
&&(\frac{\delta_{C_1}}{\sum_{\alpha}[d_{\alpha}]^2})(\frac{\delta_{C_2}}{\sum_{\alpha}[d_{\alpha}]^2})=(\frac{\delta_{C_1}}{\sum_{\alpha}[d_{\alpha}]^2}) (\frac{\delta_{C_1 \# C_2}}{\sum_{\alpha}[d_{\alpha}]^2})\nonumber
\eea
we can replace easily the correlation function by one where we put all
Lickorish generators rather than the canonical curves only. This fact will be
useful in the next section.\\

We are going to proove the last theorem through two lemmas describing some
properties of this invariant.

\bl
The expectation value
$<\prod_{i=1}^{g}\delta_{y_i}\prod_{i=1}^{g}\delta_{x_i}>_{q-YM({\cal T}_g)}$
associated to a Heegaard diagram  $({\cal T}_g, (x_j)_{j=1,...,g},
(y_j)_{j=1,...,g})$ of a manifold ${\cal M}$ is invariant under any Singer's
move applied to the diagram.
\el
\proof\\
{\bf Trivial moves:}\\
it is a fact already established that the expectation value is invariant under
any isotopic deformation of any curve in the surface, simply because of the
flatness condition.\\
We have the formula $\delta_C=\delta_{C^{-1}}$, which is exactly the second
trivial move.\\
{\bf Handlebodies Diffeomorphisms:}\\
The flatness condition implies trivially the following property for any curves
$C_1,C_2$:
\be
\delta_{C_1}\delta_{C_2}=\delta_{C_1}\delta_{C_1 \# C_2}
\ee
moreover the $\delta_{e_i}$s (resp.$\delta_{f_i}$s ) commute one with the
others. We then obtain the invariance under the first type 1 move.\\
It is easy to see, with the expression of the block element ${\cal A}_{free}$,
that we can do the following replacement along any curve $x_i$.:
\par
\centerline{\psfig{figure=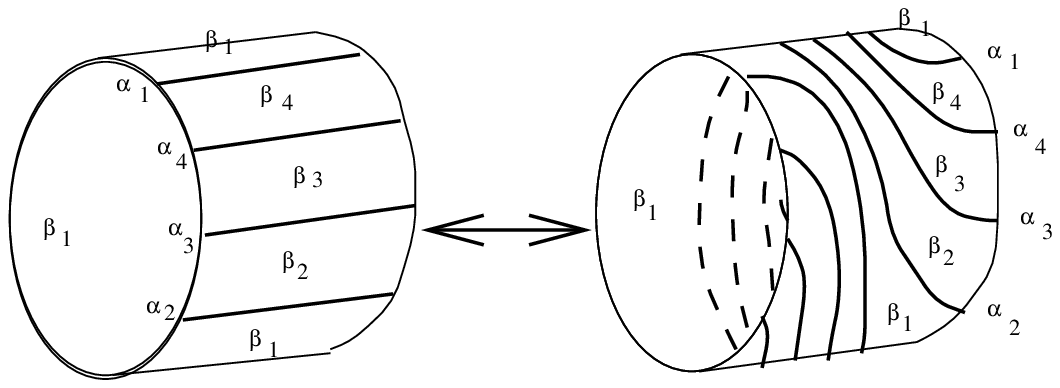}}
\par
which implies the invariance under the second type 1 move.\\
{\bf $g \rightarrow g+1$ moves:}\\
Let us consider a Heegaard diagram with one handle with its $a-$ and
$b-$cycles. We cut the surface along a certain $2-$cell to obtain a torus with
a puncture on which are drawn the two cycles as in the following figure. We
choose the minimal fat graph describing this object to describe the partial
integration of the expectation value over the edges of the latter object.
\par
\centerline{\psfig{figure=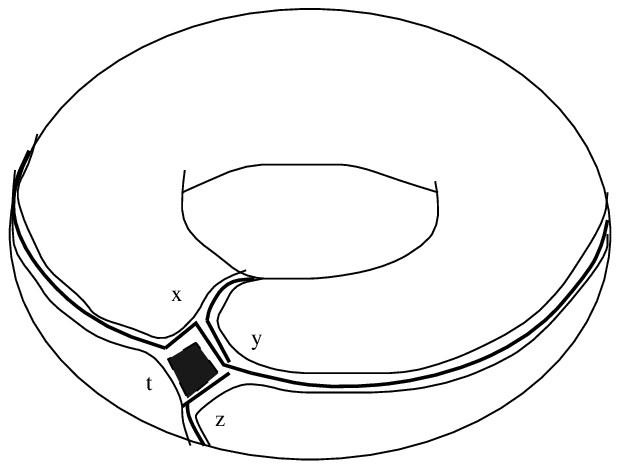}}
\par
we obtain easily:
\ben
&&\int  dh(U_{[y,t]})dh(U_{[x,z]}) \delta_{[y,t,x,z,t,y,z,x]}\delta_{[x,y,z,x]}
\delta_{[x,y,t,x]} =\\
&&= \int  dh(U_{[y,t]})dh(U_{[x,z]}) \delta_{[y,t,x,z,t]}\delta_{[x,y,z,x]}
\delta_{[x,y,t,x]}\\
&&= \int  dh(U_{[y,t]})dh(U_{[x,z]}) \delta_{[t,x,z]}\delta_{[x,y,z,x]}
\delta_{[x,y,t,x]}\\
&&= \int  dh(U_{[y,t]})dh(U_{[x,z]}) \delta_{[t,x,y,z]}\delta_{[x,y,z,x]}
\delta_{[x,y,t,x]}\\
&&=\delta_{[t,x,y,z]}\\
\een
the last line is easily obtained by using the property that the integration
just "pick" the zero component associated to a link.
The latter result establishes the invariance under the type 2 Singer move.
This ends the proof of the lemma and shows that
the expectation value is an invariant of the manifold ${\cal M}.$
\cqfd

\bl
For any augmented Heegaard diagram $({\cal L},(x_i)_{i=1 \cdots g}, (y_i)_{i=1
\cdots g}, (z_i)_{i=1 \cdots g})$ describing a manifold ${\cal M}$ there exists
a framed link $L$ which is a surgery data describing the same manifold ${\cal
M}$ and verifying:
\be
\frac{<\prod_{i=1}^{g}\delta_{y_i}\prod_{i=1}^{g}\delta_{x_i}>_{q-YM({\cal
A})}} {<\prod_{i=1}^{g}\delta_{z_i}\prod_{i=1}^{g}\delta_{x_i}>_{q-YM({\cal
A})} }=
RT({\cal M})
\ee
where RT is the Reshetikhin Turaev invariant of the manifold computed from $L.$
\el

\proof
The trick already used in the proof of the invariance under the second Singer
move can be used also here. We first use a natural property of delta functions
that can be described by the following figure :
\par
\centerline{\psfig{figure=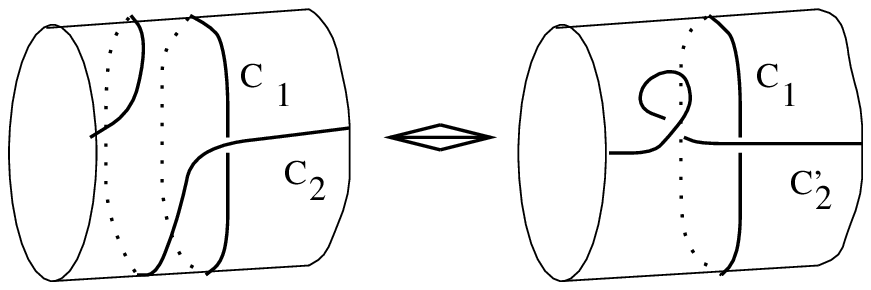}}
\par
\be
\delta_{{C}_2}\delta_{{C}_1}=\delta_{{C'}_2}\delta_{{C}_1}
\ee
to transform the correlation function in a new one\\
 $<(\sum_{\alpha_1,\cdots, \alpha_g}(\prod_{i}
[d_{\alpha_i}])W((R_i,\alpha_i)_{i=1\cdots
g})\prod_{i=1}^{g}\delta_{x_i}>_{q-YM({\cal L})}$
 where the $R_i$s are ribbons glued on the surface with the same framings and
knotted in $S^3$ in the same way as the $y_i$s but with a support now included
in the area described in the following figure:
\par
\centerline{\psfig{figure=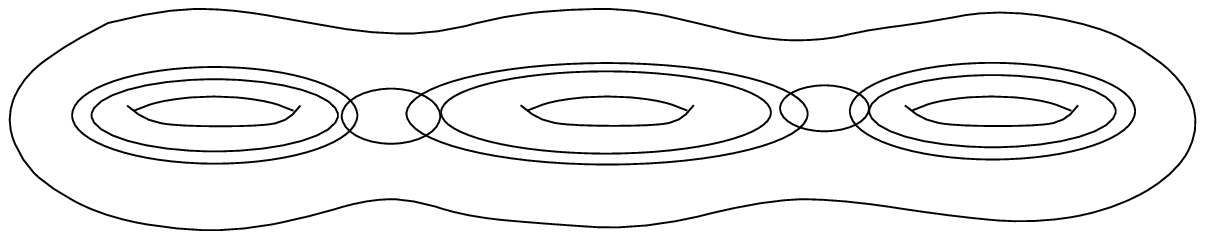}}
\par
Using now the usual flatness property (\ref{flatness}) we can deform again the
latter knot to put its crossings in the "discs", the rest of the knot being
composed of parallel strands along handles zones.\\
 Then we are able to do the same calculus as in
the verification of the invariance under type 2 Singer's move. The integration
"picks" again the zero component on each segment of the skeleton. We then
obtain the equality:
\be
<\prod_i \delta_{y_i} \prod_i \delta_{x_i}>_{q-YM({\cal
A})}=\sum_{\alpha_1,\cdots, \alpha_g}(\prod_{i} [d_{\alpha_i}])(\prod_{j}
I_{disc_j})
\ee
with $I_{disc_j}$ being the invariant associated, by our construction, to the
knot contained in the j-th disc, placed on the sphere $S^2$ and with four
coupons picking the zero component on the boundary of the disc.\\
It is  then easy to see, using the equivalence already established in section
(3) between our invariant on the sphere and the Reshetikhin invariant of link ,
that the quantity $I_{disc_j}$ is exactly the Reshetikhin invariant associated
to this framed link with coupons. \\
Now the proof can be achieved by establishing, using the Reshetikhin-Turaev
framework, that the latter data is a surgery data of the manifold ${\cal M}.$
 Let us first recall that, using the "related" surgery and Heegaard
descriptions, we can replace the set of curves $y_i$ describing the manifold
${\cal M}$ by a link composed of the curves $z_i$ associated to the Heegaard
description of $S^3$ placed at a time $t$ and the curves $R_i \times {t_i}$ (
with $t_i \le t$ ) associated to the composition of Dehn twists describing the
Heegaard gluing diffeomorphism encoded in the $y_i$s. If we compute, with the
notations of Reshetikhin and Turaev in \cite{RT2}, the invariant associated to
the framed link described before, with an insertion of two "coupons" for each
handle picking the zero component, we obtain easily that this invariant is
equal to the invariant associated to the link $L=\cup_i R_i \times {t_i}$ only.
This property uses trivially the fact that :
\ben
\sum_{\alpha,{\alpha'}_1,{\alpha'}_n}tr_{V_{\alpha}}(\mua\phi^{\alpha_n
\alpha_{n-1}}_{{\alpha'}_{n-1}} \cdots \phi^{\alpha_2
\alpha_1}_{{\alpha'}_1}\phi^{{\alpha'}_1 \alpha}_{0}
\psi_{{\alpha'}_1 \alpha}^{0}\psi_{\alpha_2
\alpha_1}^{{\alpha'}_1}\cdots\psi_{\alpha_n \alpha_{n-1}}^{{\alpha'}_{n-1}})=
id_{V_{\alpha_1}} \otimes \cdots \otimes id_{V_{\alpha_n}}
\een
We then obtain:
\ben
\frac{<\prod_{i=1}^{g}\delta_{y_i}\prod_{i=1}^{g}\delta_{x_i}>_{q-YM({\cal
L})}}{<\prod_{i=1}^{g}\delta_{z_i}\prod_{i=1}^{g}\delta_{x_i}>_{q-YM({\cal
L})}}=\sum_{(\alpha)}(\prod_{i}[d_{\alpha_i}]) RT((R_i,\alpha_i)_{i=1 \cdots
n})
\een
Using now the celebrated result of \cite{RT2} this object is a non trivial
invariant of the 3-manifold ${\cal M},$ it is the Reshetikhin-Turaev's
invariant of three manifolds. \\
\cqfd

\subsection{Chern-Simons theory on a lattice and Three dimensional Lattice
q-gauge theory }
We will define here a three dimensional gauge theory which extends in some
sense
the previous construction on a surface. The definition of this theory is based
on a choice of a simplicial presentation of the manifold which exhibits
naturally a
canonical decomposition of the manifold.
Let us consider a 3-manifold ${\cal M}$ given by a complex $K$. We impose here
that all vertices of $K$ are tetravalent. We will denote $K^{\star}$ the dual
complex of $K$.
We will denote as before $A_0^i,A_1^j,A_2^k,A_3^l$ the $0-,1-,2-,3-$simplexes
of $K$
, $A_0^{\star i},A_1^{\star j},A_2^{\star k},A_3^{\star l}$ the
$0-,1-,2-,3-$simplexes
of $K^{\star}$ and $B_1^j$ (resp.$B_1^{\star j}$) the middle of the $A_1^j$
(resp. $A_1^{\star j}$).

\bd[canonical thickening of a graph]
Let us define another tetravalent complex $ K^{\#}$ build up from the previous
one as follows:
A couple of $B_1^j$ and $B_1^{\star j}$ are said
to be a couple of neighbours if $B_1^j$ is the middle of an edge of a certain
$A_2^k$ and  $B_1^{\star j}$ is in the middle of this $A_2^k.$
We denote by $J$ the set of $1-$simplexes defined by the set of couples of
neighbour
points. We now define the $0-simplexes$ of $K^{\#}$ to be the middles of the
elements
of $J$. The $1-simplexes$ of $K^{\#}$ are then given by the set of couples
of $0-simplexes$ corresponding to elements of $J$ having one vertex in common,
if
this vertex is a $B_1^j$ (resp. a $B_1^{\star j}$) then this $1-$simplex is
said "of type $K$"
(resp. "of type $K^{\star}$").
Now the $2-$simplexes are defined to be of three types: one $2-$simplex is
associated to each closed curve formed by type $K$ $1-$simplexes only, one to
each closed
curve formed by type $K^{\star}$ $1-$simplexes only, and one to each closed
curve formed
alternatively by type $K$ and type $K^{\star}$ $1-$simplexes.We will refer us
to "the $e_{i}$s"
, "the $f_{j}$s", and "the $P$s" to denote respectively these three types of
$2-$simplexes.
Finally the $3-$simplexes are defined in an obvious way by considering each
connected region
around the vertices of $K$ and $K^{\star}$.
\ed
We will denote by $K^{\#}_0,K^{\#}_1,K^{\#}_2,K^{\#}_3$ the sets of
$0-,1-,2-,3-$simplexes respectively. A piece of this new complex is shown in
the following figure:

\par
\centerline{\psfig{figure=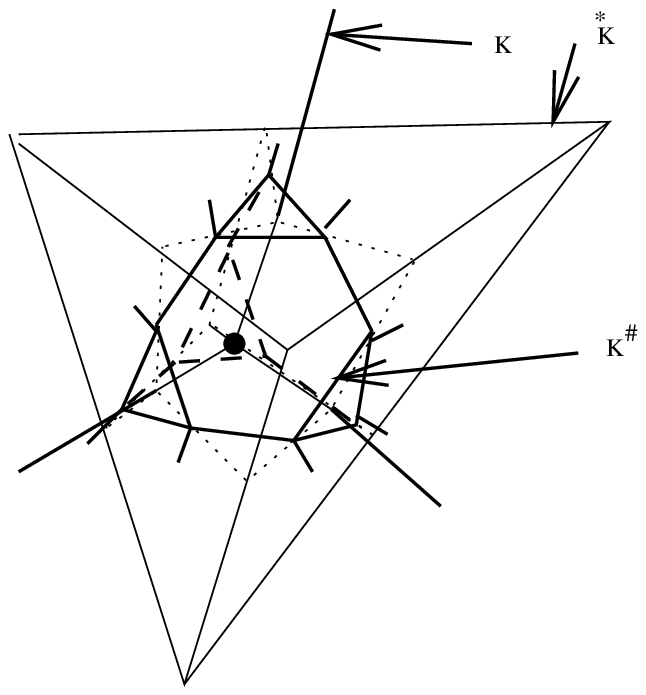}}
\par

The graph $K^{\#}$ build up from any triangulation $K$ describing a manifold
${\cal M}$ owns the
following properties:

\begin{lemma}
If we denote by ${\cal R}_1$ (resp. ${\cal R}_2$) the region defined by the set
of $3-$cells associated to vertices of $K$
(resp.$K^{\star}$) and by ${\cal L}$ the surface defined by the set of Ps.The
decomposition:
${\cal M}={\cal R}_1 + {\cal L} + {\cal R}_2$ is a canonical decomposition and
$K^{\#}$ is homeomorphic to $K$.
The set formed by the elements of $K^{\#}_0$,the elements of $ K^{\#}_1$ and
all $P$s forms the complex L
associated to the triangulation of the canonical surface {\cal L}(for this
reason these sets of $0-,1-$ and
$2-$simplexes will be also denoted respectively by $L_0,L_1,L_2$)

The set of canonical $2-$cells of ${\cal R}_1$ (resp ${\cal R}_2$) is a subset
of the $e_i$s (resp.the $f_j$s).
\end{lemma}

\proof
This decomposition is equivalent to the Heegaard decomposition
"derived" from the complex $K.$
\cqfd

\bd[3-dimensional lattice q-gauge theory]
As a consequence of the property that $K^{\#}_0=L_0$ and
$K^{\#}_1=L_1$, we can define as before the
exchange algebra associated to the elements of $K^{\#}_1$
by imposing the coaction of the gauge symmetry algebra at each element of
$K^{\#}_0$ and by choosing
a cilium order on the surface.We can define as before
the Wilson loops attached to each closed path formed by elements of $K^{\#}_1$
(i.e. drawn on L)
and delta functions associated to each $2-$cell.
In fact we define the Yang-Mills weight associated to a 2-cell $P$ of area
$A_P$ to be:
\be
\delta^{\beta}_P=\sum_{\alpha \in Phys(A)}[d_{\alpha}] e^{- \frac{A_P
C_{\alpha}}{2 \beta}}\Wa_{P}
\ee
where $C_{\alpha}$ is the quadratic casimir of the representation $\alpha$ and
$\beta$ is a coupling constant of the Yang-Mills theory.
We define the expectation value associated to any element ${\cal A}$ of
$\Lambda^{inv}$ in the 3 dimensional q-Yang Mills theory to be:
\be
<{\cal A}>_{{\cal M}}:= \int \prod_{l\in K^{\#}_1} dh(U_l) (\prod_{j}
\delta^{\beta}_{f_j})(\prod_{P\in L_2}\delta^{\beta}_P)\;\; {\cal A}\;\; (
\prod_{i} \delta^{\beta}_{e_i})
\ee
in the limit $q \rightarrow 1$ this theory becomes the well known Yang-Mills
theory on a lattice associated to a manifold ${\cal M}.$
\ed

\bp
Let $L$ be a link drawn on the 1-skeleton $K^{\#}_1$ of ${\cal M}$. Using again
the properties of the complex $K$, $L$ is in fact drawn on the canonical
surface
and we can define $W_L$ in the framework defined in this article. The
correlation function associated to $L$ in the limit $\beta \rightarrow 0$ is
then
\be
{lim}_{\beta \rightarrow \infty}\frac{<W_L>_{{\cal M}}}{<1>_{{\cal M}}} =
\frac{RT({\cal M}, L)}{RT({\cal M})}
\ee
this formula can be considered as a description of Reshetikhin-Turaev
invariants, i.e. of Chern-Simons invariants in term of a well defined lattice
gauge theory and a definition of the Witten's path integral formulas.
\ep

\proof\\
The expectation value is simply the same as that introduced in the last
subsection but with a very special Heegaard decomposition where the gluing
diffeomorphism is simply the identity. \\
\cqfd

\medskip

{\bf Acknowledgements:} It is a pleasure to thank my friend P.Roche for his
constant support. I want also aknowledge
illuminating discussions with N.Reshetikhin and C.Mercat.

\bibliographystyle{unsrt}

\end{document}